\documentclass[11pt,a4paper]{amsart}
\usepackage[T1]{fontenc}
\usepackage[utf8]{inputenc}
\usepackage{amsmath,amsfonts,amssymb,graphicx,url,amsthm,mathtools}
\numberwithin{equation}{section}
\usepackage[left=2.5cm,right=2.5cm,top=2.5cm,bottom=2.5cm]{geometry}
\usepackage{enumitem}
\usepackage{xcolor,color,nicefrac}
\usepackage{booktabs,multirow}
\usepackage{caption}
\usepackage{subcaption}
\usepackage{bm}
\usepackage[noend]{algpseudocode}
\usepackage[section,ruled]{algorithm}
\usepackage[export]{adjustbox}
\usepackage[sort]{natbib}
\usepackage{setspace}
\usepackage[hyperfootnotes=true]{hyperref}
\hypersetup{colorlinks=true,raiselinks=true,breaklinks=true,citecolor=blue}
\usepackage[capitalize,nameinlink]{cleveref}
\usepackage{tikz,pgfplots,pgfplotstable}
\pgfplotsset{compat=1.18}
\usetikzlibrary{arrows,arrows.meta,decorations,backgrounds,positioning,calc,matrix}
\usepgfplotslibrary{groupplots}
\definecolor{cSingle}{HTML}{1F77B4}   
\definecolor{cMulti}{HTML}{D62728}    
\definecolor{cGuide}{HTML}{9E9E9E}    
\definecolor{cLzero}{HTML}{2CA02C}    
\definecolor{cLa}{HTML}{1F77B4}       
\definecolor{cLb}{HTML}{FF7F0E}       
\definecolor{cLc}{HTML}{D62728}       
\definecolor{ICSgrey}{rgb}{0.3,0.3,0.35}
\pgfplotsset{
  cplot/.style={
    width=.62\linewidth, height=.45\linewidth,
    grid=both,
    minor grid style={draw=gray!12,line width=.1pt},
    major grid style={draw=gray!30,line width=.2pt},
    tick align=outside, tick pos=left,
    xlabel near ticks, ylabel near ticks,
    label style={font=\small}, tick label style={font=\footnotesize},
    legend style={font=\footnotesize,draw=gray!40,fill=white,fill opacity=.92,
                  text opacity=1,rounded corners=1pt,row sep=1pt},
    legend cell align=left,
    every axis plot/.append style={line width=0.5pt, mark size=1.5pt},
  },
}
\newtheorem{theorem}{Theorem}[section]
\newtheorem{lemma}[theorem]{Lemma}

\newtheorem{proposition}[theorem]{Proposition}

\newtheorem{remark}[theorem]{Remark}

\def\Letters{A,B,C,D,E,F,G,H,I,J,K,L,M,N,O,P,Q,R,S,T,U,V,W,X,Y,Z}
\makeatletter
\@for \@l:=\Letters \do{%
  \expandafter\edef\csname\@l bb\endcsname{\noexpand\ensuremath{%
  \noexpand\mathbb{\@l}}}%
  \expandafter\edef\csname\@l cal\endcsname{\noexpand\ensuremath{%
  \noexpand\mathcal{\@l}}}%
  \expandafter\edef\csname\@l frak\endcsname{\noexpand\ensuremath{%
  \noexpand\mathfrak{\@l}}}%
}
\makeatother
\newcommand{\isdef}{\mathrel{\mathrel{\mathop:}=}}

\newcommand{\bs}{\boldsymbol}
\newcommand{\kernel}{\mathcal{K}}
\renewcommand{\d}{\operatorname{d}}
\DeclareMathOperator{\diam}{diam}
\DeclareMathOperator{\spn}{span}
\DeclareMathOperator{\dist}{dist}
\DeclareMathOperator{\supp}{supp}
\DeclareMathOperator{\diag}{diag}
\DeclareMathOperator{\cond}{cond}
\DeclareMathOperator*{\argmin}{arg\,min}
\newcommand{\norm}[1]{\left\lVert#1\right\rVert}
\newcommand{\abs}[1]{\left\lvert#1\right\rvert}
\newcommand{\inner}[2]{\left\langle#1,#2\right\rangle}
\newcommand{\Kmat}[2]{{\bs K}_{#1,#2}}           
\newcommand{\KmatL}[2]{{\bs K}^{\Lcal}_{#1,#2}}  

\begin{document}
\title[Samplet-compressed multiscale collocation and PIKL]{Compressed Multiscale Kernel Collocation and Physics-Informed Kernel Learning for Elliptic BVPs}
\thanks{SA and MM gratefully acknowledge support from the Swiss National Science
Foundation (SNSF) through the SNSF starting grant ``Multiresolution methods for
unstructured data'' (TMSGI2\_211684).}
\author{Sara Avesani}
\address{Sara Avesani, Universit\`a della Svizzera italiana, Lugano, Switzerland}
\email{sara.avesani@usi.ch}
\author{Manel Babay}
\address{Manel Babay, Universit\`a della Svizzera italiana, Lugano, Switzerland}
\email{} 
\author{Angela Kunoth}
\address{Angela Kunoth, University of Cologne, Germany}
\email{akunoth@uni-koeln.de}
\author{Michael Multerer}
\address{Michael Multerer, Universit\`a della Svizzera italiana, Lugano, Switzerland}
\email{michael.multerer@usi.ch}
\begin{abstract}
{Kernel collocation is a meshfree method for partial differential
equations whose accuracy depends on the smoothness of the kernel and of the
solution. For globally supported kernels, however, the collocation matrices are
dense and severely ill-conditioned. We address both issues by combining a
multiscale collocation scheme with samplets, a multiresolution basis of discrete
signed measures with vanishing moments. In the samplet basis, the collocation
matrix has only $\mathcal{O}(N\log N)$ relevant entries, while the multiscale
scheme adapts the kernel width to the spacing of the points on every level. Since
the differential operator scales the equations at the interior points differently
from those at the boundary points, the rescaling alone does not bound the
condition number of the level systems. We therefore solve them by a
Schur-complement reduction that is invariant under this scaling.
For noisy data, we consider physics-informed
kernel learning, which replaces the collocation system by a regularized
least-squares fit. We show that it reduces to collocation as the regularization
vanishes, and we solve it by a deflated conjugate gradient method that reuses the
Schur-complement reduction. Numerical experiments in two and three dimensions, with exact and noisy boundary
data, show that the method handles significantly more points than dense
collocation.}
\end{abstract}
\maketitle
\section{Introduction}\label{sec:intro}

\subsection{Motivation}
We consider the numerical solution of the elliptic boundary value problem
\begin{equation}\label{eq:intro_bvp}
\begin{aligned}
\Lcal u & = f && \text{in } \Omega,\\
\Bcal u & = g && \text{on } \partial\Omega,
\end{aligned}
\end{equation}
posed on a bounded Lipschitz domain $\Omega\subset\Rbb^d$. Herein,
$\Lcal$ is a second-order
elliptic operator and $\Bcal$ is a boundary operator. Classical
discretizations such as finite elements or finite differences require a mesh,
whose generation for complicated geometries is expensive and, in higher
dimensions, becomes the dominant computational issue. This has motivated interest in \emph{meshfree} methods, where the discretization is
built directly on a set of scattered nodes without connectivity information,
see \cite{belytschko1996meshless}. Among them, strong-form collocation with
radial basis functions (RBF) is particularly attractive. {It approximates the solution by a
linear combination of RBF centered at the nodes, whose coefficients are
determined by imposing the differential equation at the interior nodes and the
boundary condition at the boundary nodes.} RBF depend only on
the Euclidean distance to a center, which makes them naturally suited to
irregular domains and scattered nodes. In addition, they provide high-order 
convergence for smooth solutions, see \cite{hardy1971multiquadric,
micchelli1986interpolation, schaback1995error, Wendland05}. The use of RBF for
partial differential equations goes back to Kansa, see
\cite{kansa1990multiquadrics, kansa1990multiquadricsI}, and the resulting
\emph{Kansa method} has since been applied to a broad range of problems, see
for instance \cite{cecil2004numerical, zerroukat1998numerical, wang2013study,
sarra2005adaptive, wang2017meshfree}.

The price of using a globally supported RBF for the discretization of
\eqref{eq:intro_bvp} is a \emph{dense} and
\emph{ill-conditioned} collocation matrix. Its condition number usually
grows rapidly when the fill distance of the nodes decreases, such that
refining the node set to improve the accuracy simultaneously degrades the
conditioning, see \cite{schaback1995error}. 
This is the \emph{trade-off} or \emph{uncertainty} principle of RBF
approximation. {Moreover, the collocation matrix is unsymmetric,
since the differential operator is applied in the rows of the interior nodes but
not in those of the boundary nodes.} The symmetric collocation scheme of
Fasshauer \cite{fasshauer1996solving}, based on Hermite interpolation, restores
symmetry but leaves both density and ill-conditioning untouched.
Together, these two
issues render kernel collocation only feasible for moderate numbers of
nodes.


\subsection{Literature review}
Several strategies have been introduced to address the dense, ill-conditioned system. One direction of work
replaces globally supported RBF by \emph{compactly supported} ones, see
\cite{schaback1993special, wendland1995piecewise}, which yields sparse
matrices at the cost of an additional support parameter that must be balanced
against accuracy. This idea is the basis of localized collocation schemes, see for instance 
\cite{divo2007efficient, kosec2008local, yao2011localized}. 

A second approach is the \emph{multiscale} method, introduced in
\cite{floater1996multistep} for interpolation, in which the residual of a
coarse approximation is successively corrected on finer node sets at smaller
RBF scales. Its convergence in Sobolev spaces on bounded domains was
established in \cite{wendland2010} for interpolation and in
\cite{farrell2013rbf} for symmetric collocation. For interpolation, scaling the
RBF proportionally to the fill distance keeps the condition number of 
the corresponding generalized Vandermonde matrix at every level uniformly bounded. 
For collocation, this does not hold. The differential operator applied to the scaled
RBF introduces an extra rescaling factor in the matrix block corresponding
to interior nodes, such that this interior block and the block of the
boundary interactions are scaled differently, see 
\cite{FarrellPestana2015,farrell2013rbf}.

In connection with the multiscale method, \emph{algebraic compression} can
be applied to the generalized Vandermonde matrices of globally supported RBF. Each of these dense matrices is compressed in a multiresolution basis, in the spirit of
wavelet compression for boundary integral operators, see
\cite{BCR, DHS, TW03}. In the scattered data setting this role is played by
\emph{samplets}, discrete signed measures with
vanishing moments, in whose basis the kernel matrix of an asymptotically
smooth RBF has only log-linearly many relevant entries and admits a
sparse matrix algebra, see \cite{HM2}. Multiscale interpolation and samplet
compression have recently been combined in \cite{avesani2024}.

The theory of unsymmetric collocation is less complete than that of the
symmetric variant. In \cite{HonSchaback2001}, it has been shown that no
general nonsingularity result can hold, by constructing point configurations
for which the Kansa matrix is singular, while observing that such configurations are rare. 
Convergence results for the unsymmetric method are available in the framework of
\cite{Schaback2007}, provided the collocation points are chosen finer than the
centers, see also \cite{liu2023unsymmetric}. As an application, the regularized
least-squares view of collocation underlies the Gaussian process solvers of
\cite{chen2021solving, owhadi2019operator} and the physics-informed learning paradigm
of \cite{raissi2019physics}, particularly in the presence of noisy data. We refer
to this formulation as \emph{physics-informed kernel learning} (PIKL).


\subsection{The approach taken}
Density and ill-conditioning in kernel collocation are typically addressed separately.
This work treats them as a single problem with a joint solution. Algebraic compression
alone is insufficient, since it replaces the exact system by a perturbed one and the solution's perturbation is amplified by the condition number. 
{The multiscale method alone improves the conditioning, but leaves
every level system dense. Moreover, in collocation, the differential operator
applied to the rescaled kernel multiplies the interior block by a factor that grows
with the level, so that the condition number of the level systems still grows.}
We therefore embed samplet compression in a multiscale collocation
scheme with the kernel scale proportional to the fill distance at each
level{, and we propose a new solver that is invariant under this
factor. At each level, the linear system is solved by a Schur-complement reduction
with respect to the interior and boundary unknowns. The factor cancels in the Schur
complement and only appears as a scaling of the interior data.}
The interior block of the collocation matrix, obtained by applying the differential
operator to the kernel, is symmetric positive definite for the Mat\'ern family. We
prove this by transporting the differentiation to the Fourier side. The interior block
is solved by the conjugate gradient method, and the small boundary Schur complement is
factorized directly. 
We show that PIKL reduces exactly to unsymmetric collocation as the
regularization parameter goes to zero, and that its regularized normal equations can
be solved by a deflated conjugate gradient method that reuses the collocation solver.


\subsection{Contributions}
The contributions of this article are the following.
\begin{itemize}[leftmargin=*]
\item We apply samplet compression to the Kansa collocation matrix and embed it
  in a multiscale collocation scheme, obtaining linear systems at each level of the
  hierarchy with
  log-linearly many nonzero entries.
\item We prove that the interior block of the collocation matrix is symmetric
  positive definite for Mat\'ern kernels of sufficient smoothness and use this
  to solve every level system by a Schur-complement reduction with a
  conjugate gradient interior solve and a direct boundary solve. The reduction
  is invariant under the scaling of the interior rows.
\item Following the analysis of symmetric collocation in \cite[Sec.~16.3]{Wendland05} and of the
Lagrange multiplier approach to append boundary conditions in \cite{dahmen2001lbb}, we separate 
the
  interior and boundary contributions to the error and derive a relation between the fill distances of the two point
  sets, together with the corresponding counting rule.
\item We show that physics-informed kernel learning with Tikhonov
  regularization reduces exactly to collocation as the regularization vanishes, and
    that under the multiscale rescaling the Tikhonov and native-space penalties
    are equivalent level by level. We solve the regularized normal equations by a
  deflated conjugate gradient method.
\item We demonstrate the method on Poisson problems in two and three spatial
  dimensions, including the volume enclosed by the Stanford bunny with many collocation points, and on noisy boundary data.
\end{itemize}

\subsection{Outline}
The remainder of this article is organized as follows. \Cref{sec:prelim}
introduces the model problem, the strong-form collocation as well as the 
preliminaries on kernel interpolation and Mat\'ern kernels that are needed later.
\Cref{sec:samplets} recalls samplets and the compression theorem.
\Cref{sec:collocation} develops the samplet-compressed collocation in its
single scale and multiscale forms, derives the balance between interior and
boundary points, and presents the Schur-complement solver.
\Cref{sec:pikl} treats physics-informed kernel learning and the solution of
its normal equations by deflation. \Cref{sec:numerics}
reports the numerical experiments, and \Cref{sec:conclusion} summarizes the
findings.

Throughout, to avoid the excessive use of generic constants, we write
$A\lesssim B$ if $A\le CB$ for a constant $C>0$ independent of the quantities
of interest, and $A\sim B$ if $A\lesssim B$ and $B\lesssim A$. The symbol
$C$ denotes a generic constant whose value may change from line to line.
{Moreover, $\norm{\cdot}_2$ denotes the Euclidean norm of vectors
and $\norm{\cdot}_F$ the Frobenius norm of matrices. The letter $N$ always
denotes the number of points of the point set under consideration.}

\section{Preliminaries}\label{sec:prelim}

\subsection{Model problem and strong-form collocation}\label{sec:bvp}
{As a model problem, we consider \eqref{eq:intro_bvp} with
$\Lcal=-\Delta$ and $\Bcal=I$, that is, the Poisson problem with Dirichlet
boundary conditions,}
\begin{equation}\label{eq:poisson}
-\Delta u = f\ \text{ in } \Omega,
\qquad
u = g\text{ on } \Gamma\isdef\partial\Omega.
\end{equation}
The weak formulation of \eqref{eq:poisson} is obtained by testing against
$v\in H^1_0(\Omega)$ and integrating by parts. Collocation instead retains the strong form and tests against Dirac
measures. Since $\delta_{\bs x}\in[H^{s}(\Omega)]'$ only for $s>d/2$, the
functional $v\mapsto(\Lcal v)(\bs x)$ is bounded on $H^\tau(\Omega)$ only if
\begin{equation}\label{eq:collocation_regularity}
\tau>2+\frac d2 .
\end{equation}
In particular $\delta_{\bs x}\notin H^{-1}(\Omega)$ for $d\ge2$. Hence, the trial
space of a collocation method must be smoother than the energy space and the
two derivatives taken before evaluation reduce the convergence rate by two orders.

Given a finite-dimensional trial space $V_N$ and a set of collocation points, we
split the collocation points into interior and boundary points,
\begin{equation}\label{eq:point_split}
X=X_I\cup X_B,\qquad X_I\subset\Omega,\qquad X_B\subset\Gamma,\qquad
N_I\isdef\abs{X_I},\qquad N_B\isdef\abs{X_B},\qquad N\isdef N_I+N_B.
\end{equation}
Strong-form collocation seeks $u_N\in V_N$ with
\begin{equation}\label{eq:collocation_equations}
(\Lcal u_N)(\bs x)=f(\bs x)\ \text{ for } \bs x\in X_I,
\qquad
(\Bcal u_N)(\bs x)=g(\bs x)\ \text{ for } \bs x\in X_B .
\end{equation}
The two sets sample objects of different dimension and regularity. How 
accurately either is resolved relative to the other is a design choice of the
method, to which we return in \Cref{sec:balance}.
{
\begin{remark}\label{rem:embedding}
Condition \eqref{eq:collocation_regularity} ensures that collocation is
well defined in the classical sense, provided that
$f\in C(\overline\Omega)$ and $g\in C(\Gamma)$. For a bounded Lipschitz domain $\Omega$ and $\tau>2+d/2$, there holds 
$H^\tau(\Omega)\hookrightarrow C^2(\overline\Omega)$, which can be shown
by standard arguments, see \cite{adams2003sobolev}.
Hence, every $v\in H^\tau(\Omega)$ and $\Lcal v$ are continuous up to
the boundary. The collocation conditions 
\eqref{eq:collocation_equations} are therefore point evaluations of continuous
functions, and the boundary trace in $\Bcal v$ reduces to the 
restriction $v|_\Gamma$.
\end{remark}
}


\subsection{Kernel interpolation and Mat\'ern kernels}\label{sec:kernels}
A symmetric function $\kernel\colon\Omega\times\Omega\to\Rbb$ is a
\emph{positive definite kernel} if, for every finite set of pairwise distinct
points $X=\{\bs x_1,\dots,\bs x_N\}\subset\Omega$, the kernel matrix
\begin{equation}\label{eq:kernel_matrix_def}
\Kmat{X}{X}\isdef\big[\kernel(\bs x_i,\bs x_j)\big]_{i,j=1}^N
\end{equation}
is symmetric positive definite. Every positive definite kernel is the
reproducing kernel of a uniquely determined Hilbert space, its \emph{native
space} $\Ncal_\kernel$, with $\inner{f}{\kernel(\cdot,\bs x)}_{\Ncal_\kernel}
=f(\bs x)$ for all $f\in\Ncal_\kernel$, see \cite{aronszajn1950theory,
Wendland05}. We consider translation invariant radial kernels
$\kernel(\bs x,\bs y)=\Phi(\bs x-\bs y)$ with $\Phi\in C(\Rbb^d)\cap
L^1(\Rbb^d)$, for which Bochner's theorem \cite{bochner1933monotone}
characterizes positive definiteness by the nonnegativity of the Fourier
transform $\widehat\Phi$. The interpolant of $u$ at $X$ is
\begin{equation}\label{eq:interpolant}
s_u\isdef\sum_{i=1}^N\alpha_i\Phi(\cdot-\bs x_i),
\qquad
\Kmat{X}{X}\bs\alpha=\big[u(\bs x_i)\big]_{i=1}^N ,
\end{equation}
and it is the norm-minimal element of $\Ncal_\kernel$ interpolating the data,
with $\norm{s_u}_{\Ncal_\kernel}^2=\bs\alpha^\top\Kmat{X}{X}\bs\alpha$.
In this context, the functions \(\Kcal(\cdot,{\bs x}_i)=\Phi(\cdot-\bs x_i)\),
\(i=1,\ldots,N\) are referred to as \emph{radial basis functions}.

The quality of the interpolant is measured against the \emph{fill distance}
and the \emph{separation radius},
\begin{equation}\label{eq:fill_distance}
h_{X,\Omega}\isdef\sup_{\bs x\in\Omega}\min_{\bs x_i\in X}\norm{\bs x-\bs x_i}_2,
\qquad
q_X\isdef\tfrac12\min_{i\ne j}\norm{\bs x_i-\bs x_j}_2 ,
\end{equation}
and $X$ is called \emph{quasi-uniform} if $h_{X,\Omega}\le c_{\mathrm{qu}}q_X$
with a constant $c_{\mathrm{qu}}$ independent of $N$. There holds the 
following error estimate from \cite[Cor.~11.33]{Wendland05}.

\begin{theorem}\label{thm:kernel_approx}
Let $\Omega\subset\Rbb^d$ be a bounded Lipschitz domain and let $\Phi$ have a
Fourier transform with algebraic decay
\begin{equation}\label{eq:matfourier}
\widehat\Phi(\bs\xi)\sim\big(1+\norm{\bs\xi}_2^2\big)^{-\tau},
\qquad\bs\xi\in\Rbb^d,\quad\tau>d/2 ,
\end{equation}
so that $\Ncal_\kernel$ is norm-equivalent to $H^\tau(\Omega)$. Then for every
$0\le t<\lceil\tau\rceil-d/2-1$, all $u\in H^\tau(\Omega)$ and all $X$ with
sufficiently small fill distance there holds
\begin{equation}\label{eq:error_bound}
\norm{u-s_u}_{H^t(\Omega)}\lesssim h_{X,\Omega}^{\tau-t}\norm{u}_{H^\tau(\Omega)} .
\end{equation}
\end{theorem}

The kernels used throughout are the \emph{Mat\'ern kernels}
\cite{MAT, Rasmussen2006}. For $\nu>0$ and a length scale $\lambda>0$,
\begin{equation}\label{eq:matern}
k_\nu(r)\isdef\frac{2^{1-\nu}}{\Gamma(\nu)}
\Bigg(\frac{\sqrt{2\nu}\,r}{\lambda}\Bigg)^{\nu}
K_\nu\Bigg(\frac{\sqrt{2\nu}\,r}{\lambda}\Bigg),
\qquad r=\norm{\bs x-\bs y}_2 ,
\end{equation}
with $K_\nu$ the modified Bessel function of the second kind. In particular, we use the Mat\'ern-5/2 kernel, which reads
\begin{equation}\label{eq:matern_halfint}
k_{5/2}(r)=\Bigg(1+\sqrt5\frac r\lambda+\frac{5r^2}{3\lambda^2}\Bigg)
\exp\Bigg(-\sqrt5\frac r\lambda\Bigg).
\end{equation}
The spectral density of $k_\nu$ is
\begin{equation}\label{eq:matern_symbol}
\widehat{k}_\nu(\bs\xi)
=\alpha_{\nu,\lambda,d}\Bigg(1+\frac{\lambda^2}{2\nu}\norm{\bs\xi}_2^2\Bigg)^{-\nu-d/2},
\end{equation}
so that \eqref{eq:matfourier} holds with $\tau=\nu+d/2$ and the native space
of $k_\nu$ is norm-equivalent to $H^{\nu+d/2}$, see \cite{Wendland05}.

Mat\'ern kernels are \emph{asymptotically smooth}, meaning that there is
$r>0$ such that
\begin{equation}\label{eq:kernel_estimate}
\biggl|\frac{\partial^{\abs{\bs\beta}+\abs{\bs\gamma}}}
{\partial\bs x^{\bs\beta}\partial\bs y^{\bs\gamma}}\Phi(\bs x-\bs y)\biggr|
\lesssim\frac{(\abs{\bs\beta}+\abs{\bs\gamma})!}
{r^{\abs{\bs\beta}+\abs{\bs\gamma}}\norm{\bs x-\bs y}_2^{\abs{\bs\beta}+\abs{\bs\gamma}}}
\end{equation}
uniformly in $\bs\beta,\bs\gamma\in\Nbb^d$ and for $\bs x\ne\bs y$. Away from 
the diagonal the kernel is thus well approximated by low-degree polynomials.

\section{Samplets and matrix compression}\label{sec:samplets}
Samplets are discrete signed measures with vanishing moments, introduced in
\cite{harbrecht2021samplets}. We recall their construction and the
compression theorem in the form needed here. The proofs can be found in
\cite{harbrecht2021samplets, HM2}.


\subsection{Samplet basis}\label{sec:samplet_basis}
Associated with the points $X=\{\bs x_1,\dots,\bs x_N\}$ are the Dirac
measures $\delta_{\bs x_i}\in[C(\Omega)]'$, which span the space
\begin{equation}\label{eq:Xspace}
\Xcal\isdef\spn\{\delta_{\bs x_1},\dots,\delta_{\bs x_N}\}
\end{equation}
of discrete signed measures supported on $X$. We equip $\Xcal$ with the inner
product defined by $\inner{\delta_{\bs x_i}}{\delta_{\bs x_j}}_\Xcal=\delta_{ij}$,
which identifies $\Xcal$ isometrically with $\Rbb^N$. A \emph{multiresolution
analysis} of $\Xcal$ is a nested sequence
$\Xcal_0\subset\Xcal_1\subset\dots\subset\Xcal_J=\Xcal$, and the increment
between two levels is the orthogonal \emph{detail space} $\Scal_j$ with
$\Xcal_{j+1}=\Xcal_j\oplus\Scal_j$. In wavelet terminology, a basis of $\Xcal_0$
consists of \emph{scaling distributions} and a basis of $\Scal_j$ of
\emph{samplets}. Together they form the \emph{samplet basis}
$\bs\Sigma=\{\sigma_{j,k}\}_{j,k}$ of $\Xcal$. A samplet has \emph{vanishing
moments} of order $q+1$ if
\begin{equation}\label{eq:vanishing_moments}
\sigma_{j,k}(p)=0\qquad\text{for all } p\in\Pcal_q ,
\end{equation}
where $\Pcal_q$ denotes the polynomials of total degree at most $q$.

The construction of samplets relies on a \emph{cluster tree} $\Tcal$ for $X$, that is, a
tree with root $X$ in which every non-leaf cluster $\tau$ is the disjoint union
of its children, obtained by recursive geometric bisection of the bounding
box. For quasi-uniform $X$ the tree is balanced with depth
$J_{\mathcal{T}}=\mathcal{O}(\log N)$. On each cluster, the samplets are obtained bottom-up
from the scaling distributions of the children by a QR factorization of the
local moment matrix of the $m_q=\binom{q+d}{d}$ monomials of degree at most
$q$. The leading columns reproduce the scaling distributions of the parent,
the remaining ones are orthonormal samplets with $q+1$ vanishing moments by
construction. This is the scattered-data analogue of the wavelet construction
of \cite{TW03}, and it costs $\mathcal{O}(Nm_q^2)$ operations. We summarize the
resulting properties, which can be found in \cite{harbrecht2021samplets}.

\begin{theorem}\label{thm:samplet_basis}
The samplet basis $\bs\Sigma$ is an orthonormal basis of $\Xcal$ such that
\begin{enumerate}[label=(\roman*)]
\item the number of samplets on level $j$ behaves like $2^j$;
\item every samplet has $q+1$ vanishing moments;
\item every samplet is supported in a single cluster $\tau$ of $\Tcal$, and for
quasi-uniform $X$ there holds $\diam\supp\sigma_{j,k}\lesssim 2^{-j/d}$.
\end{enumerate}
The change of basis from the Dirac basis to $\bs\Sigma$, the \emph{fast samplet
transform}, is orthogonal and can be applied in $\mathcal{O}(N)$ operations.
\end{theorem}

We identify the basis with the orthogonal matrix $\bs\Sigma\in\Rbb^{N\times N}$
whose columns hold the coefficients of the samplets with respect to the Dirac
basis. The vanishing moments force the samplet coefficients of a smooth
function to be small. Replacing $f$ by its Taylor polynomial on the support of
$\sigma_{j,k}$ annihilates the polynomial part and leaves a high-order
remainder. In particular, the following lemma from \cite{avesani2025multiresolution} holds.

\begin{lemma}\label{lem:decay}
Let ${\bs x}_0 \in \Omega$, $f \in C^\alpha({\bs x}_0)$ for some $\alpha \geq 0$
with $q \geq \lfloor \alpha \rfloor$. Then, for every cluster 
$\tau_{j,k}\in\Tcal$ containing ${\bs x}_0$, there holds
\begin{equation*}
|\sigma_{j,k}(f)| \lesssim \operatorname{diam}(\tau_{j,k})^\alpha
\sqrt{\#\tau_{j,k}},
\end{equation*}
where $C^\alpha({\bs x}_0)$ denotes the Jaffard regularity at ${\bs x}_0$, see \cite{jaffard1991pointwise}. 
\end{lemma}


\subsection{Samplet matrix compression}\label{sec:compression}
The kernel matrix in samplet coordinates,
\begin{equation}\label{eq:samplet_matrix}
\bs K^\Sigma\isdef\bs\Sigma^\top\Kmat{X}{X}\bs\Sigma ,
\end{equation}
is orthogonally equivalent to $\Kmat{X}{X}$ and has entries
$(\sigma_{j,k}\otimes\sigma_{j',k'})(\kernel)$. By \Cref{lem:decay} and the
asymptotic smoothness \eqref{eq:kernel_estimate}, an entry is negligible
whenever the supporting clusters are well separated relative to their sizes.
Setting such entries to zero yields the compressed matrix $\bs K^\Sigma_\rho$. The following consistency error holds, see \cite{harbrecht2021samplets}.

\begin{theorem}\label{thm:compression}
Let $X$ be quasi-uniform and set to zero every entry of $\bs K^\Sigma$ whose
supporting clusters $\tau$ and $\tau'$ satisfy the admissibility condition
\begin{equation}\label{eq:cutoff}
\dist(\tau,\tau')\ge\rho\max\{\diam(\tau),\diam(\tau')\},\qquad\rho>0 .
\end{equation}
Then the compressed matrix $\bs K^\Sigma_\rho$ satisfies
\begin{equation}\label{eq:compression_error}
\norm{\bs K^\Sigma-\bs K^\Sigma_\rho}_F
\lesssim m_q\Big(\frac{r\rho}{d}\Big)^{-2(q+1)}\norm{\bs K^\Sigma}_F ,
\end{equation}
where $r$ is the constant of \eqref{eq:kernel_estimate}, and it has
$\mathcal{O}(m_q^2N\log N)$ nonzero entries.
\end{theorem}

Note that the compression is not restricted to the
square kernel matrix. For two point sets $X,Y$ and the rectangular block
$\Kmat{X}{Y}$, the transform $\bs\Sigma_X^\top\Kmat{X}{Y}\bs\Sigma_Y$ with the
samplet bases of the row and column sets is compressed in the same way, and so
is the differentiated block
\begin{equation}\label{eq:KL_def}
\KmatL{X}{Y}\isdef\big[(\Lcal_{\bs x}\kernel)(\bs x_i,\bs y_j)\big]_{\bs x_i\in X,\,\bs y_j\in Y},
\end{equation}
since $\Lcal_{\bs x}\kernel$ consists of second derivatives of $\kernel$ and
therefore inherits the estimate \eqref{eq:kernel_estimate} with
$\abs{\bs\beta}$ shifted by two. 

\begin{figure}[htbp!]
    \centering
        \includegraphics[width=0.24\linewidth]{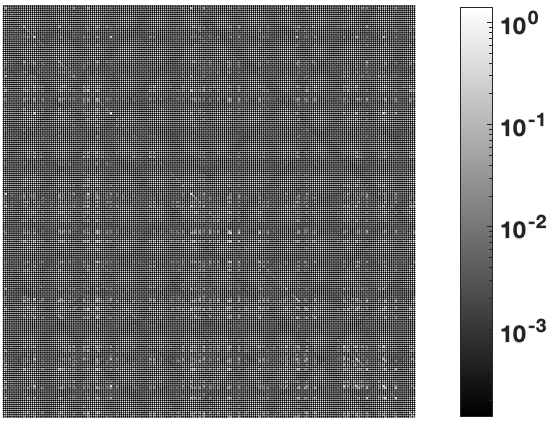} \quad \quad \quad  \quad 
       \includegraphics[width=0.24\linewidth]{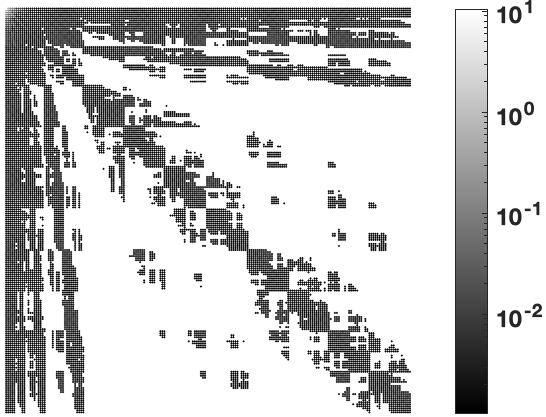}
    \caption{Visualization of the kernel matrices $\bm K^{\Sigma}$ and $\bm K^{\Sigma}_\rho $ 
    (left and right), constructed using the Matérn-$5/2$ kernel on $N = 10\,000$ points 
    uniformly distributed in the unit square $[0,1]^2$. Both matrices are displayed on 
    a logarithmic scale to reveal the large proportion of small entries, with each 
    matrix subdivided into $200 \times 200$ blocks colored by mean absolute value. 
    The transition from $\bm K^{\Sigma}$ to $\bm K^{\Sigma}_\rho $ is obtained by 
    thresholding small coefficients to zero, reducing the number of nonzero entries 
    to $\mathcal{O}(N \log N)$.}
    \label{fig:pattern_KCompressed}
\end{figure}

Moreover, for the rescaled kernel
$\Phi(\cdot/\lambda)$ the derivatives of order $\abs{\bs\beta}+\abs{\bs\gamma}$
acquire the factor $\lambda^{-\abs{\bs\beta}-\abs{\bs\gamma}}$, which cancels
against the factor $\norm{(\bs x-\bs y)/\lambda}_2^{-\abs{\bs\beta}-\abs{\bs\gamma}}$
in \eqref{eq:kernel_estimate}, so that the constant $r$ and hence the
compression rate are independent of the length scale, see \cite{avesani2024}.
This is what makes the compression applicable level by level in the multiscale
scheme of \Cref{sec:multiscale}. In practice, the compression by
\eqref{eq:cutoff} is followed by an a-posteriori thresholding of the remaining
entries at a relative tolerance $\kappa$ with respect to the Frobenius norm.
\Cref{fig:pattern_KCompressed} illustrates the sparsity pattern.

\section{Samplet-compressed kernel collocation}\label{sec:collocation}

\subsection{Kansa collocation}\label{sec:single_scale}
We solve \eqref{eq:poisson} with the
trial space spanned by the kernel translates at the collocation points,
\begin{equation}\label{eq:coll_ansatz}
u_\Phi(\bs x)=\sum_{i=1}^N\alpha_i\Phi(\bs x-\bs x_i),\qquad\bs x_i\in X ,
\end{equation}
so that the centers coincide with the collocation points and the system is
square. 
Since the Laplacian acts on the ansatz, the kernel has to satisfy
\eqref{eq:collocation_regularity}, which for Mat\'ern kernels means $\nu>2$,
see \Cref{sec:balance}. Thus, we use the Mat\'ern-$5/2$ kernel defined in
\eqref{eq:matern_halfint}. With the notation \eqref{eq:kernel_matrix_def} and
\eqref{eq:KL_def}, enforcing \eqref{eq:collocation_equations} yields the
\emph{collocation matrix}
\begin{equation}\label{eq:coll_matrix}
{\bs A}\isdef\begin{bmatrix}\KmatL{X_I}{X}\\[2pt]\Kmat{X_B}{X}\end{bmatrix}
\in\Rbb^{N\times N},
\qquad
{\bs A}\bs\alpha=\bs b\isdef\begin{bmatrix}\bs f\\ \bs g\end{bmatrix},
\end{equation}
where $\bs f$ and $\bs g$ collect $f$ at $X_I$ and $g$ at $X_B$. The matrix
$\bs A$ is unsymmetric, dense, and its condition number grows rapidly with $N$. Applying the samplet
transform of \Cref{sec:samplets} to the two blocks $\KmatL{X_I}{X}$ and
$\Kmat{X_B}{X}$ replaces $\bs A$ by a quasi-sparse matrix $\bs A_\rho$ with
$\mathcal{O}(N\log N)$ nonzero entries. Compression addresses the density but
not the conditioning. Since the compressed system is a perturbation of the
exact one, classical perturbation theory \cite{golub2013matrix} transfers the
compression error to the solution with an amplification proportional to
$\cond_2(\bs A)$, which is why the compression has to be combined with the
multiscale rescaling of \Cref{sec:multiscale}.


\subsection{Convergence and the balance of interior and boundary points}
\label{sec:balance}
For the Mat\'ern kernel $k_\nu$, whose native space is $H^{\nu+d/2}(\Omega)$,
the regularity requirement \eqref{eq:collocation_regularity} reads $\nu>2$.
This is why $k_{5/2}$ is the smallest half-integer Mat\'ern kernel admissible
for the Laplacian. Under this condition the collocation error is of the form
\begin{equation}\label{eq:coll_error}
\norm{u-u_\Phi}_{H^t(\Omega)}
\lesssim h_{X,\Omega}^{\nu+d/2-t-2}\norm{u}_{H^{\nu+d/2}(\Omega)},
\end{equation}
where the loss of two orders relative to the interpolation bound
\eqref{eq:error_bound} reflects the two derivatives taken by $\Lcal$, see
\cite{Wendland05, schaback2006kernel}.

Estimate \eqref{eq:coll_error} involves a single fill distance and does not
distinguish the two point sets of \eqref{eq:point_split}. The analysis of
\cite[Sec.~16.3]{Wendland05}, carried out for the symmetric variant of the
scheme, bounds the two components of the residual separately, each
by the fill distance of the set on which it is sampled. Let
\begin{equation}\label{eq:hI_hB}
h_I\isdef h_{X_I,\Omega},\qquad h_B\isdef h_{X_B,\Gamma}
\end{equation}
denote the fill distances of the interior points in $\Omega$ and of the
boundary points on the $(d-1)$-dimensional manifold $\Gamma$. 
In the maximum norm {on $\Omega$}, the interpolation order of
$k_\nu$ is $\nu=\tau-d/2$, which follows from \Cref{thm:kernel_approx} by
Sobolev embedding. {The same order holds at $\Gamma$. By
\Cref{rem:embedding}, we have for the Dirichlet problem that
$\Bcal(u-u_\Phi)=(u-u_\Phi)|_\Gamma$, and the
restriction to $\Gamma$ maps $H^\tau(\Omega)$ to $H^{\tau-1/2}(\Gamma)$,
reducing the rate of convergence by half an order. However, since
$\Gamma$ is a manifold of dimension $d-1$, the overall order of 
convergence is still $\tau-1/2-(d-1)/2=\nu$, see the proof of
\cite[Thm.~16.15]{Wendland05} for polygonal and \cite{FY11} for smooth
$\Gamma$.}
Since $\Lcal$ consumes two
derivatives while the Dirichlet condition consumes none, the interior residual
converges with a rate that is two orders lower than the rate of the
boundary residual, i.e.,
\begin{equation}\label{eq:wendland_split}
\norm{\Lcal(u-u_\Phi)}_{L^\infty(\Omega)}\lesssim h_I^{\nu-2}\norm{u}_{H^{\nu+d/2}(\Omega)},
\qquad
\norm{{\Bcal(u-u_\Phi)}}_{L^\infty(\Gamma)}\lesssim h_B^{\nu}\norm{u}_{H^{\nu+d/2}(\Omega)},
\end{equation}
see \cite[Thm.~16.15]{Wendland05}, and the maximum principle for elliptic
operators \cite{gilbarg2001elliptic} bounds the solution error by the sum of
the two,
\begin{equation}\label{eq:coll_error_split}
\norm{u-u_\Phi}_{L^\infty(\Omega)}\lesssim
\big(h_I^{\nu-2}+h_B^{\nu}\big)\norm{u}_{H^{\nu+d/2}(\Omega)} .
\end{equation}
For $h_I=h_B$ this recovers the rate of \eqref{eq:coll_error}, up to the order
$d/2$ that separates the maximum norm from the $L^2$-norm.

The two terms in \eqref{eq:coll_error_split} are of the same size if
\begin{equation}\label{eq:balance}
h_B\sim h_I^{(\nu-2)/\nu} .
\end{equation}
Since $(\nu-2)/\nu<1$ and $h_I<1$, the right-hand side exceeds $h_I$: at equal
error contribution the boundary is resolved more coarsely than the interior,
and the smaller $\nu$, the larger the gap. For quasi-uniform sets one has
$N_I\sim h_I^{-d}$ and $N_B\sim h_B^{-(d-1)}$, so that \eqref{eq:balance}
becomes the counting rule
\begin{equation}\label{eq:balance_counting}
N_B\sim N_I^{\frac{d-1}{d}\cdot\frac{\nu-2}{\nu}},
\end{equation}
whose exponent lies below the trace scaling $(d-1)/d$ and approaches it as
$\nu\to\infty$. A boundary discretization finer than \eqref{eq:balance} only makes 
the boundary term in \eqref{eq:coll_error_split} negligible, so that the error is
governed by the interior term alone. This is the regime in which
\eqref{eq:coll_error} applies, and it is the one used in \Cref{sec:numerics}.

\begin{remark}[Comparison with the variational setting]\label{rem:lbb}
A similar rule is known for the weak imposition of Dirichlet
conditions by Lagrange multipliers, where uniform stability of the discrete
saddle-point problem requires the trial space to be finer than the multiplier
space by a fixed number of levels $L$. For
quasi-uniform discretizations, it reads $N_B\lesssim2^{-L(d-1)}N_I^{(d-1)/d}$, see \cite{dahmen2001lbb} for details.
Note that in this variational setting a boundary finer
than the rule degrades the inf-sup constant. In collocation the boundary functionals are point evaluations,
which do not belong to $H^{-1/2}(\Gamma)$ and hence do not form a conforming
multiplier space. There is no inf-sup constant to preserve, and
\eqref{eq:balance} is a balance of accuracy.
\end{remark}

\begin{remark}[Unsymmetric collocation]\label{rem:square_case}
Estimates \eqref{eq:wendland_split}-\eqref{eq:coll_error_split} are proved
for symmetric collocation. The
matrix \eqref{eq:coll_matrix} is unsymmetric and its theory is less complete.
Invertibility cannot be guaranteed for arbitrary point configurations:
\cite{HonSchaback2001} construct singular ones, while observing that they are
rare and must be carefully constructed. The matrix \eqref{eq:coll_matrix} is
singular if and only if the trial space contains a nonzero function that is
discrete harmonic and that
vanishes at $X_B$. The convergence theory available for
unsymmetric collocation requires, as a
sufficient condition, that the collocation points be asymptotically finer than
the centers, see \cite{Schaback2007, liu2023unsymmetric} for details. Notice that this condition does not cover the square choice \eqref{eq:coll_ansatz}.
The authors of \cite{liu2023unsymmetric} note, however, that the experiments
of \cite{fasshauer2007choosing} with coinciding sets
converge although this condition is not met. In this work, we retain the square
formulation, which keeps the system of size $N$ and admits the
Schur-complement solver of \Cref{sec:schur}, and take \eqref{eq:coll_error} as
the expected rate.  We therefore assume that the Schur complement of
\Cref{sec:schur} is invertible.
\end{remark}


\subsection{Multiscale collocation}\label{sec:multiscale}
The single scale scheme uses one kernel with a fixed length scale $\lambda$
for all collocation points, and its condition number grows with $N$
\cite{kansa1990multiquadrics, fasshauer1996solving}. The multiscale approach
improves this growth by adapting the length scale level by level, see
\cite{floater1996multistep, wendland2010, farrell2013rbf, avesani2024}.
Let $J\in\Nbb$ be the number of levels. At level $\ell=1,\dots,J$ we introduce
collocation points $X_\ell=X_{\ell,I}\cup X_{\ell,B}$ with
$N_\ell\isdef\abs{X_\ell}$ and the level-dependent length scale
\begin{equation}\label{eq:multiscale_sigma}
\lambda_\ell\isdef\eta\,h_{X_\ell,\Omega},\qquad\eta>0.
\end{equation}
The rescaled kernel is then defined as
\begin{equation}\label{eq:rescaled_kernel}
\kernel_\ell(\bs x,\bs y)
\isdef \lambda_\ell^{-d} \Phi\bigg(\frac{\bs x -\bs y}{\lambda_\ell}\bigg),
\end{equation}
where $\Phi$ denotes the kernel with unit length scale.
The
blocks of the multiscale system are
\begin{equation}\label{eq:multiscale_blocks}
{\bs A}_{\ell j}\isdef
\begin{bmatrix}\KmatL{X_{\ell,I}}{X_j}\\[2pt]\Kmat{X_{\ell,B}}{X_j}\end{bmatrix},
\qquad\ell,j=1,\dots,J ,
\end{equation}
assembled from the kernel $\kernel_j$ associated with the centers of level
$j$. The multiscale approximant is the additive expansion
\begin{equation}\label{eq:multiscale_sum}
u_\Phi(\bs x)=\sum_{\ell=1}^J\sum_{i=1}^{N_\ell}
\alpha_i^{(\ell)}\kernel_\ell(\bs x,\bs x_i^{(\ell)}),
\end{equation}
whose level coefficients are computed one level at a time by solving
\begin{equation}\label{eq:multiscale_system}
{\bs A}_{\ell\ell}\bs\alpha^{(\ell)}
=\bs b_\ell-\sum_{j=1}^{\ell-1}{\bs A}_{\ell j}\bs\alpha^{(j)},
\qquad\ell=1,\dots,J ,
\end{equation}
where $\bs b_\ell$ collects the data at $X_\ell$. The right-hand side is the
residual of the data after the coarser levels have been solved. 

For interpolation, the ratio $\lambda_\ell/h_{X_\ell,\Omega}=\eta$ is constant
across levels, which bounds the condition number of every level system, see
\cite{wendland2010, farrell2013rbf}. For collocation, by the chain rule, the
Laplacian of the rescaled kernel exhibits the extra factor $\lambda_\ell^{-2}$,
which only affects the interior rows, so that the condition number of the
level systems grows with the level.


A convergence theory for the multiscale scheme was developed in
\cite{farrell2013rbf} for the symmetric variant of the level systems, building
on \cite{wendland2010, legia2012}. Provided that $u$, the coefficients of
$\Lcal$ and the boundary are sufficiently smooth, that the fill distances
decrease geometrically, $h_{X_{\ell+1},\Omega}\approx\mu\,h_{X_\ell,\Omega}$
with $\mu\in(0,1)$, and that the length scales decay slightly more slowly than
the proportional rule, $\lambda_\ell\sim h_{X_\ell,\Omega}^{1-2/(\nu+d/2)}$,
the multiscale approximant satisfies
\begin{equation}\label{eq:multiscale_error}
\norm{u-u_\Phi}_{L^2(\Omega)}\lesssim\alpha^{J}\norm{u}_{H^{\nu+d/2}(\Omega)},
\qquad\alpha\sim\mu^{\nu+d/2-2},
\end{equation}
see \cite[Thm.~4.5]{farrell2013rbf}. If $\mu$ is sufficiently small, such that
 $\alpha<1$,
the error is reduced by a constant factor with every level, and the exponent
$\nu+d/2-2$ exhibits the same loss of two orders as \eqref{eq:coll_error}. For
the unsymmetric formulation \eqref{eq:multiscale_blocks}, analogous multilevel
$L^2$ estimates have been obtained in \cite[Thm.~5.2]{liu2023unsymmetric},
under the assumption that at every level the collocation points are finer than
the centers, cf.\ \Cref{rem:square_case}.


\subsection{Schur-complement reduction}\label{sec:schur}
We now introduce the solver for the level systems, which is not affected by the different scaling of the
interior and boundary rows.

Splitting the unknowns according to $X=X_I\cup X_B$, the system takes the
block form
\begin{equation}\label{eq:schur_block}
\begin{bmatrix}\KmatL{X_I}{X_I} & \KmatL{X_I}{X_B}\\[2pt]
\Kmat{X_B}{X_I} & \Kmat{X_B}{X_B}\end{bmatrix}
\begin{bmatrix}\bs\alpha_I\\ \bs\alpha_B\end{bmatrix}
=\begin{bmatrix}\bs f\\ \bs g\end{bmatrix}.
\end{equation}
The reduction rests on the following property of the interior block.

\begin{proposition}\label{prop:derivative_pd}
Let $\Phi\in C(\Rbb^d)\cap L^1(\Rbb^d)$ be radial with
$\widehat\Phi\ge0$, $\widehat\Phi\in L^1(\Rbb^d)$ and $\widehat\Phi>0$ almost
everywhere, and let $m_1,m_2\in\Nbb_0$ with $m\isdef m_1+m_2\ge1$ satisfy the
moment condition
\begin{equation}\label{eq:pd_moment}
\int_{\Rbb^d}\norm{\bs\xi}_2^{2m}\widehat\Phi(\bs\xi)\d\bs\xi<\infty .
\end{equation}
Then $\Phi\in C^{2m}(\Rbb^d)$ and
\begin{equation}\label{eq:pd_symbol}
(-\Delta_{\bs x})^{m_1}(-\Delta_{\bs y})^{m_2}\Phi(\bs x-\bs y)=\Psi(\bs x-\bs y),
\qquad
\widehat\Psi(\bs\xi)=\norm{\bs\xi}_2^{2m}\widehat\Phi(\bs\xi) .
\end{equation}
In particular, the left-hand side depends on $m_1,m_2$ only through their sum,
and for pairwise distinct points $\bs x_1,\dots,\bs x_N$ the matrix
$[\Psi(\bs x_j-\bs x_k)]_{j,k=1}^N$ is symmetric positive definite.
\end{proposition}

\begin{proof}
Under \eqref{eq:pd_moment}, the Bochner representation
$$\Phi(\bs x)=(2\pi)^{-d/2}\int_{\Rbb^d}\widehat\Phi(\bs\xi)
e^{i\inner{\bs x}{\bs\xi}}\d\bs\xi$$ may be differentiated $2m$ times under the
integral sign, since $\norm{\bs\xi}_2^{k}\widehat\Phi\le
(1+\norm{\bs\xi}_2^2)^m\widehat\Phi\in L^1(\Rbb^d)$ for $k\le2m$ provides a
dominating function. A derivative in $x_k$ produces the factor $i\xi_k$ and a
derivative in $y_k$ the factor $-i\xi_k$; the sign is squared in a second
derivative, so that one Laplacian contributes the factor $-\norm{\bs\xi}_2^2$
in either argument. This proves \eqref{eq:pd_symbol}, and the continuity of
the dominated integrands gives $\Phi\in C^{2m}(\Rbb^d)$. Since $\widehat\Psi$
is real and even, $\Psi$ is even and the matrix is symmetric. For
$\bs c\in\Rbb^N\setminus\{\bs0\}$, exchanging the finite sum with the integral,
\[
\bs c^\top[\Psi(\bs x_j-\bs x_k)]_{j,k}\bs c
=(2\pi)^{-d/2}\int_{\Rbb^d}\norm{\bs\xi}_2^{2m}\widehat\Phi(\bs\xi)
\Big|\sum_{j=1}^N c_j e^{i\inner{\bs x_j}{\bs\xi}}\Big|^2\d\bs\xi\ge0 .
\]
As $\norm{\bs\xi}_2^{2m}\widehat\Phi>0$ almost everywhere, equality forces the
continuous function $\bs\xi\mapsto\sum_j c_j e^{i\inner{\bs x_j}{\bs\xi}}$ to
vanish identically, and the linear independence of the characters of pairwise
distinct points yields $\bs c=\bs0$.
\end{proof}

The choice $m_1=1$, $m_2=0$ is precisely the interior block of
\eqref{eq:schur_block}: for $\Lcal=-\Delta$ its symbol is
$\norm{\bs\xi}_2^2\widehat\Phi$, which gives at once the symmetry and the
positive definiteness of $\KmatL{X_I}{X_I}$. The moment condition
\eqref{eq:pd_moment} is the only restriction. For the Mat\'ern family,
whose spectral density \eqref{eq:matern_symbol} decays like
$\norm{\bs\xi}_2^{-(2\nu+d)}$, it holds if and only if $\nu>m$. 
The moment condition for a single Laplacian therefore already holds for $k_{3/2}$, whereas the doubly
differentiated kernel $\Delta_{\bs x}\Delta_{\bs y}\kernel$ of symmetric
collocation requires $\nu>2$. For the Gaussian kernel every order is allowed.

Eliminating the interior coefficients from the first block row of
\eqref{eq:schur_block} results in
\begin{equation}\label{eq:schur_zi}
\bs\alpha_I=\big(\KmatL{X_I}{X_I}\big)^{-1}\big(\bs f-\KmatL{X_I}{X_B}\bs\alpha_B\big).
\end{equation}
Then, substituting into the second row yields the Schur complement system
\begin{equation}\label{eq:schur_S}
{\bs S}\bs\alpha_B=\bs g-\Kmat{X_B}{X_I}\big(\KmatL{X_I}{X_I}\big)^{-1}\bs f,
\end{equation}
where ${\bs S}\isdef\Kmat{X_B}{X_B}-\Kmat{X_B}{X_I}\big(\KmatL{X_I}{X_I}\big)^{-1}\KmatL{X_I}{X_B}$.

The solves with the interior block are performed by the conjugate gradient
method, which is admissible by \Cref{prop:derivative_pd}.
Each matrix-vector product costs
$\mathcal{O}(N\log N)$, and the iteration count is governed by the condition
number of the interior block. The Schur complement
${\bs S}\in\Rbb^{N_B\times N_B}$ is small, since $N_B\ll N_I$, and is
factorized by a dense LU decomposition.

On level $\ell$, the blocks in \eqref{eq:schur_block} are assembled with the
rescaled kernel \eqref{eq:rescaled_kernel} on $X_I=X_{\ell,I}$ and
$X_B=X_{\ell,B}$. With $\Psi=-\Delta\Phi$, the
chain rule gives
\[
-\Delta_{\bs x}\kernel_\ell(\bs x,\bs y)
=\lambda_\ell^{-d-2}\,\Psi\Big(\frac{\bs x-\bs y}{\lambda_\ell}\Big).
\]
Hence, we can write the collocation matrix at level $\ell$ as
\begin{equation}\label{eq:level_scaling}
{\bs A}_{\ell\ell}=\bs\Theta_\ell\tilde{\bs A}_{\ell\ell},
\qquad
\bs\Theta_\ell\isdef
\begin{bmatrix}\lambda_\ell^{-2}{\bs I}_{N_I}&\bs0\\ \bs0&{\bs I}_{N_B}\end{bmatrix},
\qquad
\tilde{\bs A}_{\ell\ell}\isdef\lambda_\ell^{-d}
\begin{bmatrix}
\Psi\big(\frac{X_I-X_I}{\lambda_\ell}\big)&\Psi\big(\frac{X_I-X_B}{\lambda_\ell}\big)\\[4pt]
\Phi\big(\frac{X_B-X_I}{\lambda_\ell}\big)&\Phi\big(\frac{X_B-X_B}{\lambda_\ell}\big)
\end{bmatrix},
\end{equation}
where $\Psi\big(\frac{X-Y}{\lambda_\ell}\big)$ denotes the matrix with entries
$\Psi\big(\frac{\bs x_i-\bs y_j}{\lambda_\ell}\big)$, $\bs x_i\in X$,
$\bs y_j\in Y$, and likewise for $\Phi$. In $\tilde{\bs A}_{\ell\ell}$,
interior and boundary rows are equally rescaled. The mismatch between the two
blocks is entirely contained in $\bs\Theta_\ell$.
The next proposition quantifies the effect of $\bs\Theta_\ell$ on the
condition number of the level system.

\begin{proposition}\label{prop:schur_scaling}
Let $X_\ell$ be quasi-uniform and $\lambda_\ell\le1$. Then there is a constant
$c>0$, depending only on $\Phi$, $\eta$ and $c_{\mathrm{qu}}$, such that
\begin{equation}\label{eq:cond_level}
c\,\lambda_\ell^{-2}\le\cond_2({\bs A}_{\ell\ell})
\le\lambda_\ell^{-2}\cond_2(\tilde{\bs A}_{\ell\ell}).
\end{equation}
\end{proposition}

\begin{proof}
The upper bound follows from $\cond_2(\bs\Theta_\ell)=\lambda_\ell^{-2}$. For
the lower bound, every row $\bs r$ of ${\bs A}_{\ell\ell}$ satisfies
$\sigma_{\min}({\bs A}_{\ell\ell})\le\norm{\bs r}_2\le\sigma_{\max}({\bs A}_{\ell\ell})$.
An interior row contains the diagonal entry $\lambda_\ell^{-d-2}\Psi(\bs0)>0$,
whereas a boundary row has norm at most $C\lambda_\ell^{-d}$, with $C$
independent of $\ell$. Hence
$\cond_2({\bs A}_{\ell\ell})\ge\Psi(\bs0)\lambda_\ell^{-2}/C$.
\end{proof}

If $\cond_2(\tilde{\bs A}_{\ell\ell})$ remains bounded across levels, the
condition number of the level system grows like
$\lambda_\ell^{-2}=(\eta h_{X_\ell,\Omega})^{-2}$, due to
$\bs\Theta_\ell$. 
The Schur-complement reduction overcomes this limitation. Denoting by
$\tilde{\bs A}_{II}$, $\tilde{\bs A}_{IB}$, $\tilde{\bs A}_{BI}$,
$\tilde{\bs A}_{BB}$ the blocks of $\tilde{\bs A}_{\ell\ell}$ and inserting
\eqref{eq:level_scaling} into \eqref{eq:schur_zi} and \eqref{eq:schur_S} gives
\begin{equation}\label{eq:schur_scaled}
\bs\alpha_I=\tilde{\bs A}_{II}^{-1}\big(\lambda_\ell^{2}\bs f
-\tilde{\bs A}_{IB}\bs\alpha_B\big),
\qquad
{\bs S}\bs\alpha_B=\bs g-\tilde{\bs A}_{BI}\tilde{\bs A}_{II}^{-1}\lambda_\ell^{2}\bs f,
\qquad
{\bs S}=\tilde{\bs A}_{BB}-\tilde{\bs A}_{BI}\tilde{\bs A}_{II}^{-1}\tilde{\bs A}_{IB}.
\end{equation}
The factor $\lambda_\ell^{-2}$ cancels in the Schur complement and appears
only as the factor $\lambda_\ell^{2}$ on the interior data.

\section{Physics-informed kernel learning}\label{sec:pikl}
Collocation enforces the differential equation and the boundary condition
exactly at the collocation points. Physics-informed kernel learning replaces
these constraints by a weighted least-squares problem with an additional
regularization term. We first formulate this problem and relate it to
collocation. We then show that it can be solved with the Schur-complement
reduction of \Cref{sec:schur}.

The coefficients $\hat{\bs\alpha}\in\Rbb^N$ of \eqref{eq:coll_ansatz} are
determined by
\begin{equation}\label{eq:pikl}
\hat{\bs \alpha}
=
\argmin_{\bs \alpha\in\Rbb^N}
\left[
\frac{1}{N_I}
\left\|
\KmatL{X_I}{X}\bs \alpha-\bs f
\right\|_2^2
+
\frac{\Lambda_1}{N_B}
\left\|
\Kmat{X_B}{X}\bs \alpha-\bs g
\right\|_2^2
+
\Lambda_2\norm{{\bs R}\bs \alpha}_2^2
\right],
\end{equation}
where $\Lambda_1>0$ balances the interior and boundary residuals and
$\Lambda_2\geq0$ weights the regularization. The choice
${\bs R}^\top{\bs R}={\bs K}$ penalizes the native-space norm, since
$\bs\alpha^\top{\bs K}\bs\alpha=\norm{u_\Phi}_{\Ncal_\kernel}^2$. We use
${\bs R}={\bs I}$ instead, that is, Tikhonov regularization of the coefficient
vector. 
Indeed, under the multiscale rescaling \eqref{eq:rescaled_kernel}, both
choices are equivalent on every level: for quasi-uniform points, the kernel
matrix of level $\ell$ satisfies
\[
c_1\lambda_\ell^{-d}\norm{\bs\alpha}_2^2\le\bs\alpha^\top{\bs K}\bs\alpha
\le c_2\lambda_\ell^{-d}\norm{\bs\alpha}_2^2,
\]
with constants depending only on the kernel, on $\eta$ and on the
quasi-uniformity constant, see \cite[Chap.~12]{Wendland05} and
\cite[Sec.~3]{avesani2024}. The difference between the two regularizers is
thus absorbed into $\Lambda_2$ on each level. 

With the collocation matrix $\bs A$ and the data vector $\bs b$ of
\eqref{eq:coll_matrix}, and the weight matrix
\begin{equation}\label{eq:weight_matrix}
    \bs W\isdef\diag\Big(N_I^{-1}\bs I_{N_I},\frac{\Lambda_1}{N_B}\bs I_{N_B}\Big), 
\end{equation} 
using ${\bs R}={\bs I}$, problem \eqref{eq:pikl} reads
\begin{equation}\label{eq:link-J}
\hat{\bs\alpha}
=
\argmin_{\bs\alpha\in\Rbb^N}
\left[
(\bs A\bs\alpha-\bs b)^\top
{\bs W}
(\bs A\bs\alpha-\bs b)
+
\Lambda_2\norm{\bs\alpha}_2^2
\right],
\end{equation}
and its unique minimizer solves the normal equations
\begin{equation}\label{eq:link-normal}
\left(
{\bs A}^\top{\bs W}{\bs A}
+
\Lambda_2{\bs I}
\right)
\hat{\bs\alpha}
=
{\bs A}^\top{\bs W}\bs b .
\end{equation}
If the regularization vanishes, PIKL reduces to collocation. In particular, the following proposition holds.

\begin{proposition}\label{prop:pikl-to-colloc}
If ${\bs A}$ is nonsingular, then, for every $\Lambda_1>0$,
\[
\lim_{\Lambda_2\to0^+}\hat{\bs\alpha}
=
\bs\alpha_0
\isdef
{\bs A}^{-1}\bs b.
\]
\end{proposition}

\begin{proof}
Since ${\bs A}$ is nonsingular and ${\bs W}$ is positive definite,
${\bs A}^\top{\bs W}{\bs A}$ is nonsingular, and letting $\Lambda_2\to0^+$ in
\eqref{eq:link-normal} yields
$({\bs A}^\top{\bs W}{\bs A})^{-1}{\bs A}^\top{\bs W}\bs b={\bs A}^{-1}\bs b$.
\end{proof}

Moreover, for small $\Lambda_2>0$, the regularization perturbs the collocation
coefficients to first order. Expanding \eqref{eq:link-normal} in $\Lambda_2$
gives
\begin{equation}\label{eq:pikl_perturbation}
\hat{\bs\alpha}-\bs\alpha_0
=
-\Lambda_2
({\bs A}^\top{\bs W}{\bs A})^{-1}
\bs\alpha_0
+
\mathcal{O}(\Lambda_2^2).
\end{equation}

\subsection{Solution of the regularized system}\label{sec:pikl_solver}
Conjugate gradient applied directly to \eqref{eq:link-normal} would require
products with ${\bs A}^\top{\bs W}{\bs A}$, whose conditioning
depends on the diagonal in \({\bs W}\). We
therefore solve for the weighted collocation values
$\bs y\isdef{\bs W}^{1/2}{\bs A}\bs\alpha$. Multiplying \eqref{eq:link-normal}
by ${\bs W}^{-1/2}{\bs A}^{-\top}$ gives
\begin{equation}\label{eq:pikl_y}
\Big({\bs I}+\Lambda_2{\bs W}^{-1/2}({\bs A}{\bs A}^\top)^{-1}{\bs W}^{-1/2}\Big)\bs y
={\bs W}^{1/2}\bs b,
\qquad
\bs\alpha={\bs A}^{-1}{\bs W}^{-1/2}\bs y .
\end{equation}
The matrix in \eqref{eq:pikl_y} is symmetric positive definite. Moreover, it
is never assembled, since every product with
$({\bs A}{\bs A}^\top)^{-1}={\bs A}^{-\top}{\bs A}^{-1}$ amounts to two
Schur-complement reductions of \Cref{sec:schur}.

To solve \eqref{eq:pikl_y} efficiently, we split $\bs y=[\bs y_I;\bs y_B]$ into interior and boundary unknowns. By
\eqref{eq:level_scaling}, the interior block of $({\bs A}{\bs A}^\top)^{-1}$
contains the factor $\lambda_\ell^{4}$. As a consequence, the interior diagonal block of
the matrix in \eqref{eq:pikl_y} is close to the identity on fine levels. Its
boundary diagonal block, instead, grows with $N_B/\Lambda_1$, so that the
eigenvalues that increase the iteration count of conjugate gradients are
associated with the boundary unknowns. 
We remove them by deflation
\cite{nicolaides1987deflation, saad2000deflated, tang2009comparison}. The
block
$\bs y_B$ is eliminated exactly, and conjugate gradients are applied only to
the resulting Schur complement for $\bs y_I$. The boundary diagonal block is a small
dense $N_B\times N_B$ matrix, which is assembled once from the
Schur-complement reduction and factorized directly. \Cref{alg:pikl} summarizes the main steps of the method.

\begin{algorithm}[htb]
\caption{Deflated conjugate gradient method for \eqref{eq:link-normal}.}
\label{alg:pikl}
\begin{algorithmic}[1]
\Require compressed blocks of ${\bs A}$, weights ${\bs W}$, data $\bs b$,
parameter $\Lambda_2$, tolerance
\State factorize ${\bs A}$ by the Schur-complement reduction of
\Cref{sec:schur}
\State assemble the boundary diagonal block of \eqref{eq:pikl_y} and compute
its Cholesky factorization
\State eliminate $\bs y_B$ from \eqref{eq:pikl_y}
\State solve the resulting system for $\bs y_I$ by conjugate gradients, where
each product with the matrix costs two solves with ${\bs A}$ and ${\bs A}^\top$
and one with the boundary diagonal block
\State recover $\bs y_B$ from the boundary rows of \eqref{eq:pikl_y}
\State $\bs\alpha\gets{\bs A}^{-1}{\bs W}^{-1/2}\bs y$
\end{algorithmic}
\end{algorithm}

\section{Numerical results}\label{sec:numerics}
This section validates the theory on Poisson problems in two and three
spatial dimensions. All computations were performed on a MacBook~Pro with an
Apple~M2~Max processor and 32\,GB of RAM, using the \texttt{C++} library
\texttt{FMCA} \cite{fmca}.

\subsection{Setting and parameters}\label{sec:parameters}
The boundary value problems are solved with the multiscale, samplet-compressed
collocation scheme of \Cref{sec:multiscale}, using the Schur-complement
reduction of \Cref{sec:schur} with a conjugate gradient (CG) interior solve,
and with the physics-informed kernel learning of \Cref{sec:pikl}, using
\Cref{alg:pikl}. The kernel is the Mat\'ern-$5/2$ kernel
\eqref{eq:matern_halfint}, the least smooth half-integer Mat\'ern kernel
admissible for a second-order operator. The samplet basis is built with vanishing moments of order $q+1$, with $q=3$
in two and $q=2$ in three dimensions. The
admissibility parameter of \eqref{eq:cutoff} is $\rho=0.5$ in two and
$\rho=0.7$ in three dimensions, followed by an a-posteriori thresholding of
the compressed blocks at the relative tolerance $\kappa=10^{-6}$ and
$\kappa=10^{-5}$ with respect to the Frobenius norm, respectively. A nugget of $10^{-10}$ is added to the interior block, and the multiscale
length scales follow \eqref{eq:multiscale_sigma} with $\eta=3$ in two and
$\eta=1.5$ in three dimensions. The tolerances of the conjugate gradient
iterations are stated with each experiment.

We quantify the accuracy by the relative $\ell^2$ error on a fixed evaluation
set $X_{\mathrm{eval}}$,
\begin{equation}\label{eq:eps2}
\varepsilon_2\isdef
\frac{\norm{u-u_\Phi}_{\ell^2(X_{\mathrm{eval}})}}{\norm{u}_{\ell^2(X_{\mathrm{eval}})}},
\end{equation}
where $u_\Phi$ is the computed approximant, together with
the observed order of convergence
\begin{equation}\label{eq:order}
\mathrm{order}_\ell\isdef
\log\Bigg(\frac{\varepsilon_2^{(\ell)}}{\varepsilon_2^{(\ell-1)}}\Bigg)\Big/
\log\Bigg(\frac{h_{X_\ell,\Omega}}{h_{X_{\ell-1},\Omega}}\Bigg).
\end{equation}
{The error $\varepsilon_2$ is a discrete counterpart of the
relative $L^2(\Omega)$ error, so that the observed orders are compared with
the $L^2$ estimates \eqref{eq:coll_error}, for $t=0$, and
\eqref{eq:multiscale_error}.}

\subsection{Two-dimensional Poisson problems}\label{sec:numerics_2d}
We consider two Poisson problems on the unit square $\Omega=[0,1]^2$, chosen
to test the method in two regularity regimes. The first is
\begin{equation}\label{num:problem1}
-\Delta u_1(x,y)=2\big(x(1-x)+y(1-y)\big)\ \text{ in }\Omega,
\qquad u_1=0\ \text{ on }\partial\Omega,
\end{equation}
with the smooth polynomial solution $u_1(x,y)=x(1-x)y(1-y)$, which lies in
$H^\tau(\Omega)$ for every $\tau$. The second is designed to limit the regularity of the solution,
\begin{equation}\label{num:problem2}
-\Delta u_2(x,y)=-\tfrac{15}{4}\abs{x-0.6}^{1/2}\ \text{ in }\Omega,
\qquad u_2=\abs{x-0.6}^{5/2}\ \text{ on }\partial\Omega,
\end{equation}
with exact solution $u_2(x,y)=\abs{x-0.6}^{5/2}$, whose third derivative is
singular along the line $x=0.6$, so that $u_2\in H^s(\Omega)$ only for $s<3$.

The aim is to exhibit the stability of the multiscale system and the
scalability of the solver by comparison with classical single scale
collocation. We test \emph{single scale} (SS), consisting of a
single point cloud with the fixed length scale $\lambda_{\mathrm{ss}}=0.3$,
versus \emph{multiscale} (MS), with point sets $X_1,\dots,X_J$ of decreasing
fill distance and $\lambda_\ell=\eta h_{X_\ell,\Omega}$. We also compare the \emph{dense} ($\delta$) versus
\emph{compressed} ($\Sigma$) system. 
The point sets are uniform grids of $N_{\ell,I}=m^2$ interior points with
$N_{\ell,B}=m$ equally spaced boundary points, so that the regime of \Cref{sec:balance} applies. The levels are listed in
\Cref{tab:levels2d}. The evaluation set {$X_{\mathrm{eval}}$ is a uniform
$1000\times1000$ grid in $\Omega$}. 
The interior block is solved by CG to a relative
residual of $10^{-6}$, with the iteration count capped at $1\,000$. The dense
assembly stores the full $N\times N$ blocks and becomes infeasible beyond
$N\approx50\,000$, while the compressed scheme reaches $N\approx500\,000$ on
the same machine. To compare the schemes at equal computational effort we
also report $\varepsilon_2$ against the cost proxy
$\mathrm{cost}=N\cdot n_{\mathrm{CG}}$, where $n_{\mathrm{CG}}$ is the number of conjugate gradient iterations.

\begin{table}[htb]
\centering\small
\begin{tabular}{cccccc}
\toprule
$\ell$ & $N_{\ell,I}$ & $N_{\ell,B}$ & $N_\ell$ & $h_{X_\ell,\Omega}$ & $\lambda_\ell$ \\
\midrule
1 &     400 &  20 &     420 & $4.76\cdot10^{-2}$ & $1.43\cdot10^{-1}$ \\
2 &   1\,600 &  40 &   1\,640 & $2.44\cdot10^{-2}$ & $7.32\cdot10^{-2}$ \\
3 &   6\,400 &  80 &   6\,480 & $1.23\cdot10^{-2}$ & $3.70\cdot10^{-2}$ \\
4 &  25\,600 & 160 &  25\,760 & $6.21\cdot10^{-3}$ & $1.86\cdot10^{-2}$ \\
5 & 102\,400 & 320 & 102\,720 & $3.12\cdot10^{-3}$ & $9.35\cdot10^{-3}$ \\
6 & 409\,600 & 640 & 410\,240 & $1.56\cdot10^{-3}$ & $4.68\cdot10^{-3}$ \\
\bottomrule
\end{tabular}
\caption{Numbers of interior, boundary and total points, fill distance and
kernel length scale per level for the two-dimensional experiments.}
\label{tab:levels2d}
\end{table}

\subsubsection{Smooth solution}
\begin{figure}[htbp]
\centering
\begin{subfigure}[t]{0.48\textwidth}
\centering
\begin{tikzpicture}
\begin{loglogaxis}[cplot,
  width=\linewidth, height=5.5cm,
  xlabel={$N$}, ylabel={$\varepsilon_2$},
  ymin=1e-5, ymax=1e1,
  legend style={font=\small, draw=none, fill=none,
                at={(0.70,0.35)}, anchor=south west, cells={anchor=west}},
]
\addplot[cSingle,solid,mark=*]
coordinates{(420,1.384e-2)(1640,3.063e-3)(6480,4.151e-4)(25760,2.525e-4)};
\addlegendentry{SS$_\delta$}
\addplot[cSingle,densely dashed,mark=o,mark options={solid}]
  coordinates{(420,1.384e-2)(1640,3.056e-3)(6480,4.798e-4)(25760,8.664e-1)(102720,2.350e0)};
\addlegendentry{SS$_\Sigma$}
\addplot[cMulti,solid,mark=square*]
coordinates{(420,5.788e-3)(1640,2.462e-3)(6480,8.678e-4)(25760,2.865e-4)};
\addlegendentry{MS$_\delta$}
\addplot[cMulti,densely dashed,mark=diamond,mark options={solid,line width=.8pt}]
coordinates{(420,5.788e-3)(1640,2.462e-3)(6480,8.676e-4)(25760,2.867e-4)(102720,9.076e-5)(410240,2.765e-5)};
\addlegendentry{MS$_\Sigma$}
\draw[gray, densely dotted, line width=0.8pt]
  (axis cs:350,3e-2) -- (axis cs:600000,1.13e-4)
  node[pos=0.1, sloped, above, font=\tiny] {$\mathcal{O}(N^{-3/4})$};
\end{loglogaxis}
\end{tikzpicture}
\caption{}
\end{subfigure}
\hfill
\begin{subfigure}[t]{0.48\textwidth}
\centering
\begin{tikzpicture}
\begin{loglogaxis}[cplot,
  width=\linewidth, height=5.5cm,
  xlabel={cost $=N\cdot n_{\mathrm{CG}}$}, ylabel={$\varepsilon_2$},
  ymin=1e-5, ymax=1e1,
  legend style={font=\small, draw=none, fill=none,
                at={(0.70,0.6)}, anchor=north west, cells={anchor=west}},
]
\addplot[cSingle,solid,mark=o,mark options={solid}]
  coordinates{(30660,1.384e-2)(328000,3.056e-3)(2643840,4.798e-4)(25760000,8.664e-1)(102720000,2.350e0)};
\addlegendentry{SS$_\Sigma$}
\addplot[cMulti,densely dashed,mark=diamond,mark options={solid,line width=.8pt}]
coordinates{(20160,5.788e-3)(121360,2.462e-3)(518400,8.676e-4)(2138080,2.867e-4)(8628480,9.076e-5)(34870400,2.765e-5)};
\addlegendentry{MS$_\Sigma$}
\end{loglogaxis}
\end{tikzpicture}
\caption{}
\end{subfigure}
\caption{Relative $\ell^2$ error as a function of the problem size (a) and of the
solver work (b) for the Poisson problem \eqref{num:problem1}. The legend reads as
follows: SS$_\delta$ for the single scale dense, SS$_\Sigma$ for the single scale
samplet-compressed, MS$_\delta$ for the multiscale dense, and MS$_\Sigma$ for the
multiscale samplet-compressed scheme.}
\label{fig:coll_errors_prob1}
\end{figure}
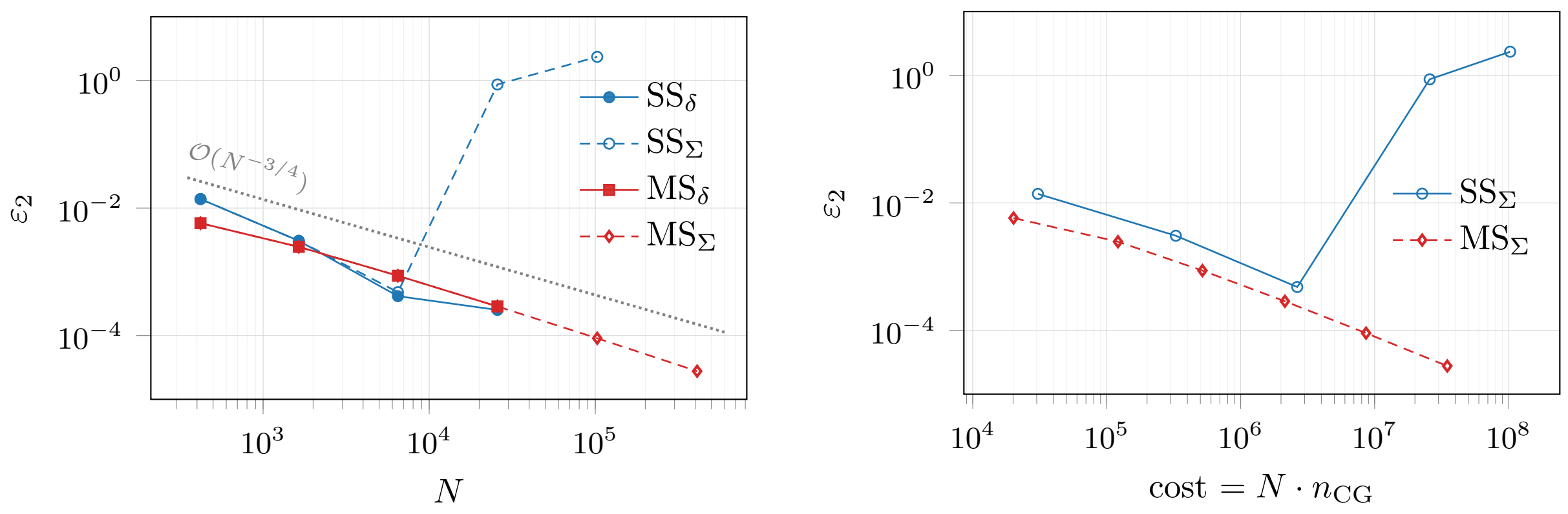

The results for \eqref{num:problem1} are collected in \Cref{tab:res_prob1} and
shown in \Cref{fig:coll_errors_prob1}. The multiscale dense and compressed
errors coincide over the whole range in which the dense assembly is feasible,
meaning that the compression preserves the accuracy of the exact scheme while
reducing the storage to $\mathcal{O}(N\log N)$.

The theory guarantees the rate $h^{\tau-2}=h^{3/2}$, with $\tau=\nu+d/2$, that is
$\mathcal{O}(N^{-3/4})$ in two dimensions, both for the single scale estimate
\eqref{eq:coll_error} and, level by level, for the multiscale one
\eqref{eq:multiscale_error}. The multiscale scheme meets this rate from the
third level on and exceeds it on the finer levels. The single scale scheme
initially produces higher orders, since the fixed width
$\lambda_{\mathrm{ss}}=0.3$ is large compared with the fill distance, and a
wide Mat\'ern kernel approaches spectral accuracy on a smooth solution. The
order drops at $N=25\,760$, which is an effect of the trade-off principle.

The same mechanism explains the failure of the compressed single scale scheme,
which follows the dense one up to $N=1\,640$, deviates from it at $N=6\,480$
and diverges from $N=25\,760$ on. The perturbation introduced by the
compression is amplified in the solution proportionally to the condition
number of the collocation matrix. The multiscale scheme avoids this
amplification, since the interior block of the Schur-complement reduction of
\Cref{sec:schur} is the rescaled matrix $\tilde{\bs A}_{II}$, which depends
only on the rescaled points and is equally well conditioned on every level.
This effect is equally visible in the iteration counts, which reach the cap of
$1\,000$ from $N=25\,760$ on for the single scale scheme, while the iteration
count of the multiscale solver levels off over three orders of magnitude
in $N$.

\begin{table}[htb]
\centering\footnotesize
\begin{tabular}{rr rrrr rrrr}
\toprule
 & & \multicolumn{4}{c}{single scale} & \multicolumn{4}{c}{multiscale}\\
\cmidrule(lr){3-6}\cmidrule(lr){7-10}
$\ell$ & $N_\ell$
 & $\varepsilon_2^{\delta}$ & $\varepsilon_2^{\Sigma}$ & order & CG
 & $\varepsilon_2^{\delta}$ & $\varepsilon_2^{\Sigma}$ & order & CG\\
\midrule
1 &      420 & $1.38\cdot10^{-2}$ & $1.38\cdot10^{-2}$ & --   &   73
            & $5.79\cdot10^{-3}$ & $5.79\cdot10^{-3}$ & --   &   48\\
2 &   1\,640 & $3.06\cdot10^{-3}$ & $3.06\cdot10^{-3}$ & 2.25 &  200
            & $2.46\cdot10^{-3}$ & $2.46\cdot10^{-3}$ & 1.28 &   74\\
3 &   6\,480 & $4.15\cdot10^{-4}$ & $4.80\cdot10^{-4}$ & 2.94 &  408
            & $8.68\cdot10^{-4}$ & $8.68\cdot10^{-4}$ & 1.53 &   80\\
4 &  25\,760 & $2.52\cdot10^{-4}$ & $8.66\cdot10^{-1}$ & 0.72 & 1000
            & $2.87\cdot10^{-4}$ & $2.87\cdot10^{-4}$ & 1.61 &   83\\
5 & 102\,720 & --                 & $2.35\cdot10^{0}$  & --   & 1000
            & --                 & $9.08\cdot10^{-5}$ & 1.67 &   84\\
6 & 410\,240 & --                 & --                 & --   & --
            & --                 & $2.77\cdot10^{-5}$ & 1.72 &   85\\
\bottomrule
\end{tabular}
\caption{Results for the Poisson problem \eqref{num:problem1}: relative $\ell^2$
error \eqref{eq:eps2} for the dense ($\delta$) and samplet-compressed ($\Sigma$)
variants, observed order \eqref{eq:order} and number of CG iterations of the
interior solve of the compressed variant. A dash marks the resolutions at which
the scheme is no longer feasible.}
\label{tab:res_prob1}
\end{table}

\subsubsection{Solution of reduced regularity}
\begin{figure}[htbp]
\centering
\begin{subfigure}[t]{0.48\textwidth}
\centering
\begin{tikzpicture}
\begin{loglogaxis}[cplot,
  width=\linewidth, height=5.5cm,
  xlabel={$N$}, ylabel={$\varepsilon_2$},
  ymin=1e-4, ymax=1e2,
  legend style={font=\small, draw=none, fill=none,
                at={(0.70,0.3)}, anchor=south west, cells={anchor=west}},
]
\addplot[cSingle,solid,mark=*]
coordinates{(420,6.258e-2)(1640,9.950e-3)(6480,1.227e-3)(25760,9.205e-4)};
\addlegendentry{SS$_\delta$}
\addplot[cSingle,densely dashed,mark=o,mark options={solid}]
coordinates{(420,6.258e-2)(1640,9.960e-3)(6480,1.616e-3)(25760,3.766e-1)(102720,1.738e1)};
\addlegendentry{SS$_\Sigma$}
\addplot[cMulti,solid,mark=square*]
coordinates{(420,1.250e-1)(1640,4.617e-2)(6480,1.603e-2)(25760,5.240e-3)};
\addlegendentry{MS$_\delta$}
\addplot[cMulti,densely dashed,mark=diamond,mark options={solid,line width=.8pt}]
coordinates{(420,1.250e-1)(1640,4.617e-2)(6480,1.603e-2)(25760,5.240e-3)(102720,1.647e-3)(410240,4.824e-4)};
\addlegendentry{MS$_\Sigma$}
\draw[gray, densely dotted, line width=0.8pt]
  (axis cs:350,5e-1) -- (axis cs:600000,1.88e-3)
  node[pos=0.1, sloped, above, font=\tiny] {$\mathcal{O}(N^{-3/4})$};
\end{loglogaxis}
\end{tikzpicture}
\caption{}
\end{subfigure}
\hfill
\begin{subfigure}[t]{0.48\textwidth}
\centering
\begin{tikzpicture}
\begin{loglogaxis}[cplot,
  width=\linewidth, height=5.5cm,
  xlabel={cost $=N\cdot n_{\mathrm{CG}}$}, ylabel={$\varepsilon_2$},
  ymin=1e-4, ymax=1e2,
  legend style={font=\small, draw=none, fill=none,
                at={(0.40,0.9)}, anchor=north west, cells={anchor=west}},
]
\addplot[cSingle,solid,mark=o,mark options={solid}]
  coordinates{(38640,6.258e-2)(510040,9.960e-3)(6480000,1.616e-3)(25760000,3.766e-1)(102720000,1.738e1)};
\addlegendentry{SS$_\Sigma$}
\addplot[cMulti,densely dashed,mark=diamond,mark options={solid,line width=.8pt}]
coordinates{(21840,1.250e-1)(123000,4.617e-2)(505440,1.603e-2)(2086560,5.240e-3)(8525760,1.647e-3)(33639680,4.824e-4)};
\addlegendentry{MS$_\Sigma$}
\end{loglogaxis}
\end{tikzpicture}
\caption{}
\end{subfigure}
\caption{Relative $\ell^2$ error as a function of the problem size (a) and of the
solver work (b) for the Poisson problem \eqref{num:problem2}. The legend reads as
follows: SS$_\delta$ for the single scale dense, SS$_\Sigma$ for the single scale
samplet-compressed, MS$_\delta$ for the multiscale dense, and MS$_\Sigma$ for the
multiscale samplet-compressed scheme.}
\label{fig:coll_errors_prob2}
\end{figure}
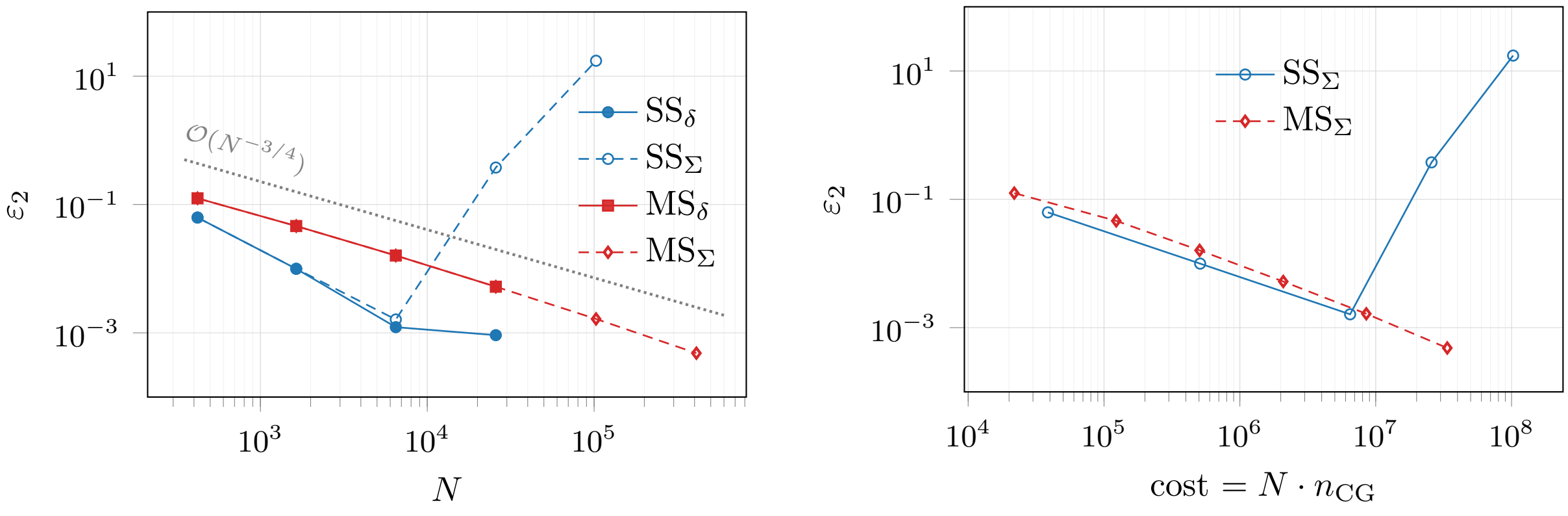

The results for \eqref{num:problem2} are collected in \Cref{tab:res_prob2} and
shown in \Cref{fig:coll_errors_prob2}. The error level is more than one order
of magnitude larger than for \eqref{num:problem1} at every $N$, so that the
singularity costs accuracy uniformly. The observed multiscale orders are nevertheless
comparable to those of the smooth problem. This is due to the fact that the
singularity is confined to a line, which contributes little to the averaged
error \eqref{eq:eps2}.

The main difference from the smooth case lies in the single scale scheme,
which fails here in both its variants. The dense variant gains high orders at
the first two refinements, followed by a substantial drop, at which point the
error stagnates near $10^{-3}$. The compressed variant follows it up to
$N=1\,640$, deviates from it at $N=6\,480$, where the iteration count reaches
the cap, and diverges from $N=25\,760$ on. The multiscale scheme, by contrast,
converges monotonically over the entire range, with an iteration count that
levels off as for the smooth problem.

\begin{table}[htb]
\centering\footnotesize
\begin{tabular}{rr rrrr rrrr}
\toprule
 & & \multicolumn{4}{c}{single scale} & \multicolumn{4}{c}{multiscale}\\
\cmidrule(lr){3-6}\cmidrule(lr){7-10}
$\ell$ & $N_\ell$
 & $\varepsilon_2^{\delta}$ & $\varepsilon_2^{\Sigma}$ & order & CG
 & $\varepsilon_2^{\delta}$ & $\varepsilon_2^{\Sigma}$ & order & CG\\
\midrule
1 &      420 & $6.26\cdot10^{-2}$ & $6.26\cdot10^{-2}$ & --   &   92
            & $1.25\cdot10^{-1}$ & $1.25\cdot10^{-1}$ & --   &   52\\
2 &   1\,640 & $9.95\cdot10^{-3}$ & $9.96\cdot10^{-3}$ & 2.75 &  311
            & $4.62\cdot10^{-2}$ & $4.62\cdot10^{-2}$ & 1.49 &   75\\
3 &   6\,480 & $1.23\cdot10^{-3}$ & $1.62\cdot10^{-3}$ & 3.07 & 1000
            & $1.60\cdot10^{-2}$ & $1.60\cdot10^{-2}$ & 1.55 &   78\\
4 &  25\,760 & $9.21\cdot10^{-4}$ & $3.77\cdot10^{-1}$ & 0.42 & 1000
            & $5.24\cdot10^{-3}$ & $5.24\cdot10^{-3}$ & 1.63 &   81\\
5 & 102\,720 & --                 & $1.74\cdot10^{1}$  & --   & 1000
            & --                 & $1.65\cdot10^{-3}$ & 1.68 &   83\\
6 & 410\,240 & --                 & --                 & --   & --
            & --                 & $4.82\cdot10^{-4}$ & 1.78 &   82\\
\bottomrule
\end{tabular}
\caption{Results for the Poisson problem \eqref{num:problem2}, at the same
resolutions and in the same format as \Cref{tab:res_prob1}.}
\label{tab:res_prob2}
\end{table}

\subsection{Three-dimensional problem on the Stanford bunny}\label{sec:bunny}
We turn to a three-dimensional domain, considering the volume
$\Omega$ enclosed by the Stanford bunny of the \texttt{Stanford 3D Scanning
Repository}, with $\partial\Omega$ its surface. We solve $-\Delta u=f$ in
$\Omega$ with $u=g$ on $\partial\Omega$, using as exact solution the
superposition of two Gaussian sources placed inside the body,
\begin{equation}\label{eq:bunny_sol}
u(\bs x)=\sum_{k=1}^2\exp\Bigg(-\frac{\norm{\bs x-\bs p_k}_2^2}{2s_k^2}\Bigg),
\qquad
\Delta u(\bs x)=\sum_{k=1}^2\exp\Bigg(-\frac{\norm{\bs x-\bs p_k}_2^2}{2s_k^2}\Bigg)
\Bigg(\frac{\norm{\bs x-\bs p_k}_2^2}{s_k^4}-\frac3{s_k^2}\Bigg),
\end{equation}
with $\bs p_1=(-1.51,2.07,0.40)$, $s_1=0.6$ in the head and
$\bs p_2=(0,0.5,0)$, $s_2=0.4$ in the body, and $g=u|_{\partial\Omega}$.
Interior and boundary points are extracted from a volumetric and a surface
sampling of the bunny, comprising $674\,787$ and $139\,262$ points,
respectively, and the error is evaluated on the union of the two clouds,
comprising $814\,049$ points.

\begin{figure}[htbp]
\centering
\includegraphics[width=0.16\linewidth]{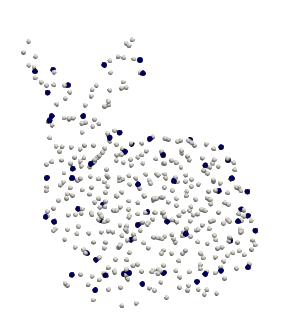}\quad\quad\quad
\includegraphics[width=0.16\linewidth]{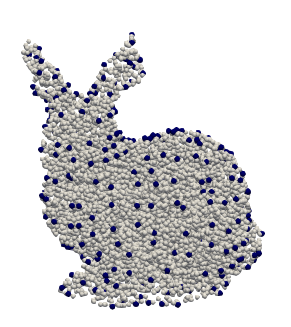}\quad\quad\quad
\includegraphics[width=0.16\linewidth]{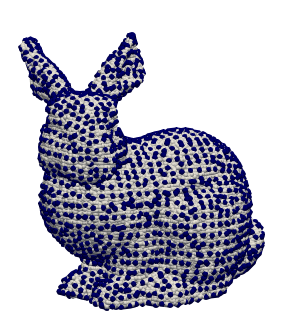}
\caption{Point sets of the hierarchy on the Stanford bunny at the
levels $\ell=2$, $\ell=4$ and $\ell=6$ (from left to right), with $454$,
$6\,745$ and $104\,589$ points.}
\label{fig:bunny_levels}
\end{figure}

The hierarchy is built by farthest-point sampling: each cloud is ordered so
that the next point is the one farthest from all points chosen so far, and
level $\ell$ consists of the first $N_{\ell,I}$ interior and $N_{\ell,B}$
boundary points of the respective ordering, with
$N_{\ell,B}\sim N_{\ell,I}^{2/3}$. {This boundary sampling is
finer than the counting rule \eqref{eq:balance_counting} requires, so that
the boundary term in \eqref{eq:coll_error_split} is negligible, see
\Cref{sec:balance}.}
The levels are nested by construction, and the separation
radius of every level is at least half its fill distance, so that the mesh
ratio is at most $2$. The fill distance $h_\ell$ is the covering radius of the
level with respect to the full cloud, and the length scale is
$\lambda_\ell=\eta h_\ell$ with $\eta=1.5$, see \Cref{tab:bunny_levels} and
\Cref{fig:bunny_levels}. The interior block is solved by CG to a relative
residual of $10^{-5}$, with the iteration count capped at $200$.

\begin{table}[htb]
\centering\small
\begin{tabular}{ccccc}
\toprule
$\ell$ & $N_{\ell,I}$ & $N_{\ell,B}$ & $h_\ell$ & $\lambda_\ell$ \\
\midrule
1 &      100 &     30 & $4.06\cdot10^{-1}$ & $6.09\cdot10^{-1}$ \\
2 &      400 &     54 & $2.35\cdot10^{-1}$ & $3.53\cdot10^{-1}$ \\
3 &   1\,600 &    137 & $1.40\cdot10^{-1}$ & $2.10\cdot10^{-1}$ \\
4 &   6\,400 &    345 & $8.36\cdot10^{-2}$ & $1.25\cdot10^{-1}$ \\
5 &  25\,600 &    869 & $5.00\cdot10^{-2}$ & $7.50\cdot10^{-2}$ \\
6 & 102\,400 & 2\,189 & $2.92\cdot10^{-2}$ & $4.38\cdot10^{-2}$ \\
\bottomrule
\end{tabular}
\caption{Nested farthest-point hierarchy on the Stanford bunny: interior and
boundary points, fill distance and length scale per level.}
\label{tab:bunny_levels}
\end{table}

\begin{figure}[htbp]
\centering
\includegraphics[width=0.24\linewidth]{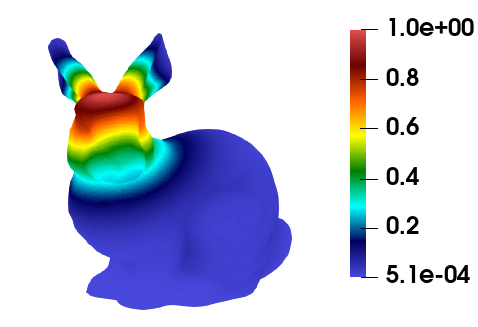}\hfill
\includegraphics[width=0.24\linewidth]{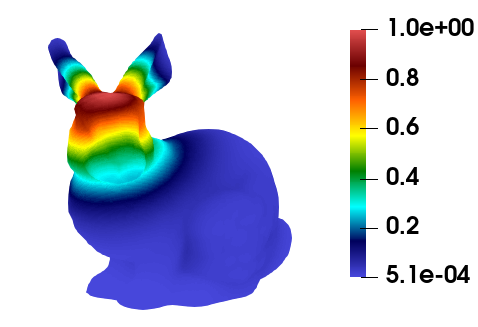}\hfill
\includegraphics[width=0.24\linewidth]{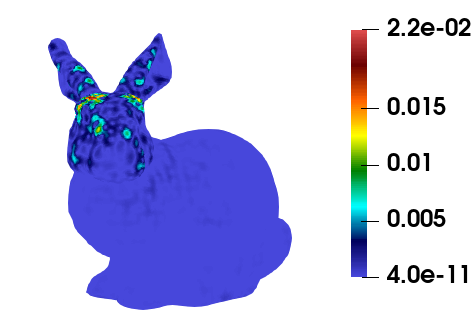}\hfill
\includegraphics[width=0.24\linewidth]{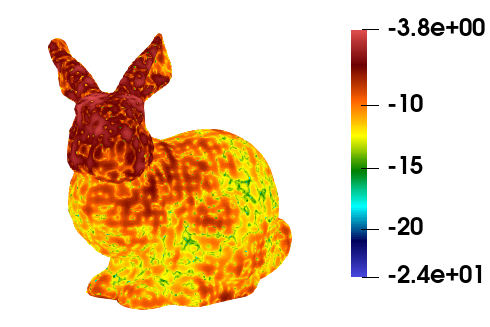} \\
\includegraphics[width=0.24\linewidth]{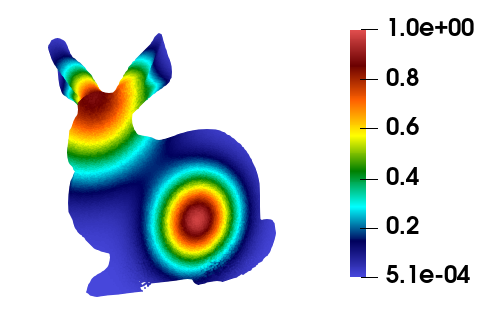}\hfill
\includegraphics[width=0.24\linewidth]{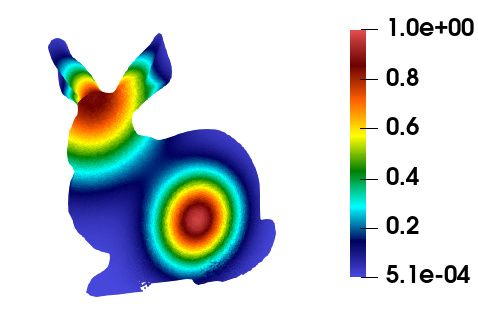}\hfill
\includegraphics[width=0.24\linewidth]{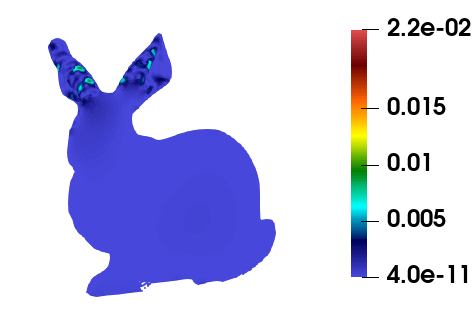}\hfill
\includegraphics[width=0.24\linewidth]{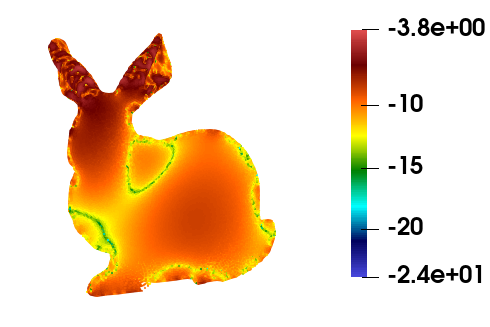}
\caption{Stanford bunny at the finest level $\ell=6$. From left to right:
exact solution \eqref{eq:bunny_sol}, multiscale samplet-compressed
collocation approximant, absolute pointwise error $\abs{u-u_\Phi}$, and the same error on a
logarithmic scale. The bottom row shows the clip with center
$(-0.3,1.2,0.1)^\top$ in direction $(0,0,1)^\top$ of the first row.}
\label{fig:bunny_solution}
\end{figure}

\Cref{tab:bunny} collects the results, and \Cref{fig:bunny_solution} shows the
solution and the error at the finest level. The observed orders are below the
rate $h^{\tau-2}=h^2$ predicted for the Mat\'ern-$5/2$ kernel in three
dimensions. This is due to the small ratio $\eta$. For multiscale
interpolation, the error is reduced from one level to the next by a factor of
the form $C(\mu^{\tau}+\eta^{-\tau})$, see \cite{wendland2010}: the first term
gives the rate in $h$, whereas the second depends only on $\eta$ and dominates
when $\eta\mu<1$. With the collocation order $\tau-2$ in place of $\tau$, the
second term equals $\eta^{-2}\approx0.44$ on the bunny, where
$\eta\mu\approx0.9$, and agrees with the observed reduction of the error per
level. In two dimensions, $\eta\mu=1.5$, and the predicted rate is attained.
The number of CG iterations stays bounded over the whole hierarchy, since the
rescaled interior block $\tilde{\bs A}_{II}$ depends only on the rescaled
points, whose mesh ratio is bounded.

\begin{table}[htb]
\centering\small
\begin{tabular}{rrccrr}
\toprule
$\ell$ & $N_{\ell}$ & $\varepsilon_2$ & order & \%nz & CG \\
\midrule
1 &      130 & $5.43\cdot10^{-2}$ & --   & 100.0 & 13 \\
2 &      454 & $2.57\cdot10^{-2}$ & 1.37 &  98.7 & 17 \\
3 &   1\,737 & $9.98\cdot10^{-3}$ & 1.83 &  81.9 & 19 \\
4 &   6\,745 & $4.59\cdot10^{-3}$ & 1.50 &  43.9 & 20 \\
5 &  26\,469 & $1.91\cdot10^{-3}$ & 1.70 &  17.8 & 19 \\
6 & 104\,589 & $8.59\cdot10^{-4}$ & 1.49 &   4.7 & 19 \\
\bottomrule
\end{tabular}
\caption{Multiscale samplet-compressed collocation on the Stanford bunny:
relative $\ell^2$ error \eqref{eq:eps2}, observed order \eqref{eq:order},
percentage of nonzero entries of the compressed level matrix and number of CG
iterations of the interior solve.}
\label{tab:bunny}
\end{table}

\subsection{Physics-informed kernel learning}\label{sec:numerics_pikl}
We now test the regularized formulation of \Cref{sec:pikl} with
${\bs R}={\bs I}$, solved by \Cref{alg:pikl}, with the weight matrix defined in~\eqref{eq:weight_matrix}
and $\Lambda_1=10^2$. The first
experiment, on exact data, verifies that PIKL reduces to collocation as
$\Lambda_2\to0$, as stated in \Cref{prop:pikl-to-colloc}, and measures what
is lost in accuracy and gained in iterations for $\Lambda_2>0$. The second
experiment is performed on noisy boundary data.

\subsubsection{Exact data}
We solve \eqref{num:problem1} with the hierarchy of \Cref{tab:levels2d} up to
$N=102\,720$. All conjugate gradient iterations, the interior ones of the
Schur-complement reductions and the one of \Cref{alg:pikl}, are run to a
relative residual of $10^{-8}$. The error is evaluated on a uniform
$1000\times1000$ grid, and the times are cumulative over the levels.

\Cref{tab:pikl_equiv} compares collocation with PIKL at $\Lambda_2=0$. Both
runs use the same point sets, the same compressed blocks and the same
evaluation grid, so that they differ only in the solver. The errors agree to
four significant digits at every level. This confirms \Cref{prop:pikl-to-colloc} numerically. 
The conjugate gradient method of \Cref{alg:pikl} terminates after two steps at
every level. For $\Lambda_2=0$, the matrix in \eqref{eq:pikl_y} is the
identity, and one step would suffice in exact arithmetic. In our
implementation, however, every step of this outer iteration applies the
Schur-complement reduction, which solves with the interior block by an inner
conjugate gradient iteration. Since the inner and the outer iterations stop at
the same tolerance, the identity is reproduced only up to this tolerance, and
a second step is needed.
PIKL is only slightly slower than collocation, since the compression of the
blocks and the Schur-complement reduction, shared by both methods, dominate
the cost of a level.

\begin{table}[htb]
\centering\footnotesize
\begin{tabular}{rr ccc cc cc}
\toprule
 & & \multicolumn{3}{c}{$\varepsilon_2$}
   & \multicolumn{2}{c}{iterations}
   & \multicolumn{2}{c}{time [s]}\\
\cmidrule(lr){3-5}\cmidrule(lr){6-7}\cmidrule(lr){8-9}
$\ell$ & $N_\ell$ & collocation & PIKL & rel.\ diff.
 & collocation & PIKL & collocation & PIKL\\
\midrule
1 &      420 & $5.787830\cdot10^{-3}$ & $5.787828\cdot10^{-3}$ & $3.5\cdot10^{-7}$ &  69 & 2 &   0.8 &   0.9\\
2 &   1\,640 & $2.461892\cdot10^{-3}$ & $2.461838\cdot10^{-3}$ & $2.2\cdot10^{-5}$ & 103 & 2 &   2.9 &   3.1\\
3 &   6\,480 & $8.676329\cdot10^{-4}$ & $8.676453\cdot10^{-4}$ & $1.4\cdot10^{-5}$ & 111 & 2 &  13.0 &  14.0\\
4 &  25\,760 & $2.866568\cdot10^{-4}$ & $2.866631\cdot10^{-4}$ & $2.2\cdot10^{-5}$ & 116 & 2 &  85.1 &  89.8\\
5 & 102\,720 & $9.076273\cdot10^{-5}$ & $9.076494\cdot10^{-5}$ & $2.4\cdot10^{-5}$ & 118 & 2 & 801.2 & 827.2\\
\bottomrule
\end{tabular}
\caption{Collocation against PIKL with $\Lambda_2=0$ on the Poisson problem
\eqref{num:problem1}: relative $\ell^2$ error \eqref{eq:eps2}, interior CG
iterations of collocation and CG iterations of \Cref{alg:pikl}, and
cumulative time.}
\label{tab:pikl_equiv}
\end{table}

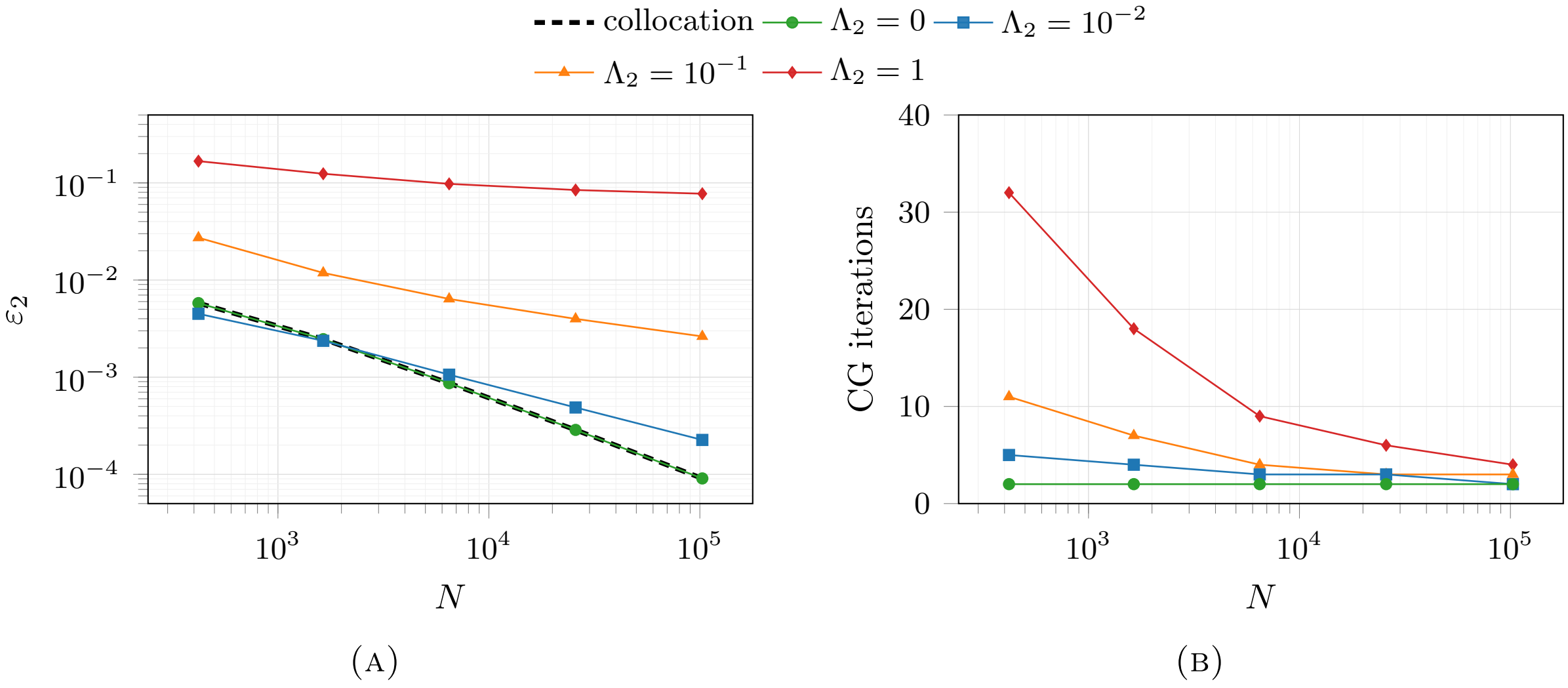
\begin{figure}[htbp]
\centering
\begin{subfigure}[t]{0.48\textwidth}
\centering
\begin{tikzpicture}
\begin{loglogaxis}[cplot,
  width=\linewidth, height=5.5cm,
  xlabel={$N$}, ylabel={$\varepsilon_2$},
  ymin=5e-5, ymax=5e-1,
  legend style={font=\footnotesize, draw=none, fill=none,
                at={(1.15,1.03)}, anchor=south, legend columns=3},
]
\addplot[black,densely dashed,line width=1.4pt,mark=none]
  coordinates{(420,5.787830e-3)(1640,2.461892e-3)(6480,8.676329e-4)(25760,2.866568e-4)(102720,9.076273e-5)};
\addlegendentry{collocation}
\addplot[cLzero,solid,mark=*]
coordinates{(420,5.787828e-3)(1640,2.461838e-3)(6480,8.676453e-4)(25760,2.866631e-4)(102720,9.076494e-5)};
\addlegendentry{$\Lambda_2=0$}
\addplot[cLa,solid,mark=square*]
  coordinates{(420,4.494858e-3)(1640,2.368898e-3)(6480,1.061168e-3)(25760,4.875615e-4)(102720,2.260228e-4)};
\addlegendentry{$\Lambda_2=10^{-2}$}
\addplot[cLb,solid,mark=triangle*]
  coordinates{(420,2.726318e-2)(1640,1.185020e-2)(6480,6.389977e-3)(25760,3.983229e-3)(102720,2.628400e-3)};
\addlegendentry{$\Lambda_2=10^{-1}$}
\addplot[cLc,solid,mark=diamond*]
  coordinates{(420,1.674257e-1)(1640,1.241467e-1)(6480,9.767784e-2)(25760,8.446347e-2)(102720,7.736578e-2)};
\addlegendentry{$\Lambda_2=1$}
\end{loglogaxis}
\end{tikzpicture}
\caption{}
\end{subfigure}
\hfill
\begin{subfigure}[t]{0.48\textwidth}
\centering
\begin{tikzpicture}
\begin{semilogxaxis}[cplot,
  width=\linewidth, height=5.5cm,
  xlabel={$N$}, ylabel={CG iterations},
  ymin=0, ymax=40,
]
\addplot[cLc,solid,mark=diamond*]
  coordinates{(420,32)(1640,18)(6480,9)(25760,6)(102720,4)};
\addplot[cLb,solid,mark=triangle*]
  coordinates{(420,11)(1640,7)(6480,4)(25760,3)(102720,3)};
\addplot[cLa,solid,mark=square*]
  coordinates{(420,5)(1640,4)(6480,3)(25760,3)(102720,2)};
\addplot[cLzero,solid,mark=*]
  coordinates{(420,2)(1640,2)(6480,2)(25760,2)(102720,2)};
\end{semilogxaxis}
\end{tikzpicture}
\caption{}
\end{subfigure}
\caption{PIKL on the Poisson problem \eqref{num:problem1} for four values of the
regularization weight: relative $\ell^2$ error against the number of
points (a) and number of CG iterations of \Cref{alg:pikl} (b).}
\label{fig:pikl_lambda}
\end{figure}

\Cref{tab:pikl_lambda} and \Cref{fig:pikl_lambda} report the same hierarchy
for $\Lambda_2\in\{0,10^{-2},10^{-1},1\}$. Since the data are exact, every
$\Lambda_2>0$ introduces a bias with respect to collocation, quantified by
\eqref{eq:pikl_perturbation}. From $N=6\,480$ on, the error indeed grows
monotonically with $\Lambda_2$. The iteration counts grow
with $\Lambda_2$ but decrease with the level: the boundary unknowns are
eliminated exactly, and the remaining interior block of \eqref{eq:pikl_y}
approaches the identity as the level increases, see \Cref{sec:pikl_solver}.

\begin{table}[htb]
\centering\footnotesize
\begin{tabular}{rr cc cc cc cc}
\toprule
 & & \multicolumn{2}{c}{$\Lambda_2=0$} & \multicolumn{2}{c}{$\Lambda_2=10^{-2}$}
   & \multicolumn{2}{c}{$\Lambda_2=10^{-1}$} & \multicolumn{2}{c}{$\Lambda_2=1$}\\
\cmidrule(lr){3-4}\cmidrule(lr){5-6}\cmidrule(lr){7-8}\cmidrule(lr){9-10}
$\ell$ & $N_\ell$
 & $\varepsilon_2$ & it & $\varepsilon_2$ & it
 & $\varepsilon_2$ & it & $\varepsilon_2$ & it\\
\midrule
1 &     420 & $5.79\cdot10^{-3}$ & 2 & $4.49\cdot10^{-3}$ &  5
            & $2.73\cdot10^{-2}$ & 11 & $1.67\cdot10^{-1}$ & 32\\
2 &   1\,640 & $2.46\cdot10^{-3}$ & 2 & $2.37\cdot10^{-3}$ &  4
            & $1.19\cdot10^{-2}$ &  7 & $1.24\cdot10^{-1}$ & 18\\
3 &   6\,480 & $8.68\cdot10^{-4}$ & 2 & $1.06\cdot10^{-3}$ &  3
            & $6.39\cdot10^{-3}$ &  4 & $9.77\cdot10^{-2}$ &  9\\
4 &  25\,760 & $2.87\cdot10^{-4}$ & 2 & $4.88\cdot10^{-4}$ &  3
            & $3.98\cdot10^{-3}$ &  3 & $8.45\cdot10^{-2}$ &  6\\
5 & 102\,720 & $9.08\cdot10^{-5}$ & 2 & $2.26\cdot10^{-4}$ &  2
            & $2.63\cdot10^{-3}$ &  3 & $7.74\cdot10^{-2}$ &  4\\
\bottomrule
\end{tabular}
\caption{PIKL on the Poisson problem \eqref{num:problem1} for four values of the
regularization weight: relative $\ell^2$ error \eqref{eq:eps2} and number of
CG iterations of \Cref{alg:pikl}.}
\label{tab:pikl_lambda}
\end{table}

\subsubsection{Noisy boundary data on the Stanford bunny}
We finally consider the geometry and the hierarchy of \Cref{sec:bunny}, 
where the boundary datum is replaced by $g(\bs x)+e(\bs x)$ with independent
$e(\bs x)\sim\mathcal{N}(0,\sigma_e^2)$. The standard deviation
\[
\sigma_e=0.05\,\frac{\norm{g}_{\ell^2(\partial\Omega)}}{\sqrt{M}}\approx0.013
\]
is $5\%$ of the root mean square of the exact datum over the full surface
cloud of $M=139\,262$ points. 
One noise realization is drawn on the
full cloud and restricted to the levels, so that a point has the same noisy
value on every level on which it appears. All conjugate gradient iterations
are run to a relative residual of $10^{-6}$. We compare
$\Lambda_2\in\{0,10^{-2},1\}$, where $\Lambda_2=0$ is collocation on the
noisy data.

\begin{figure}[htbp]
\centering
\begin{tabular}{@{}ccc@{}}
$\Lambda_2=0$ & $\Lambda_2=10^{-2}$ & $\Lambda_2=1$\\
\includegraphics[width=0.3\linewidth]{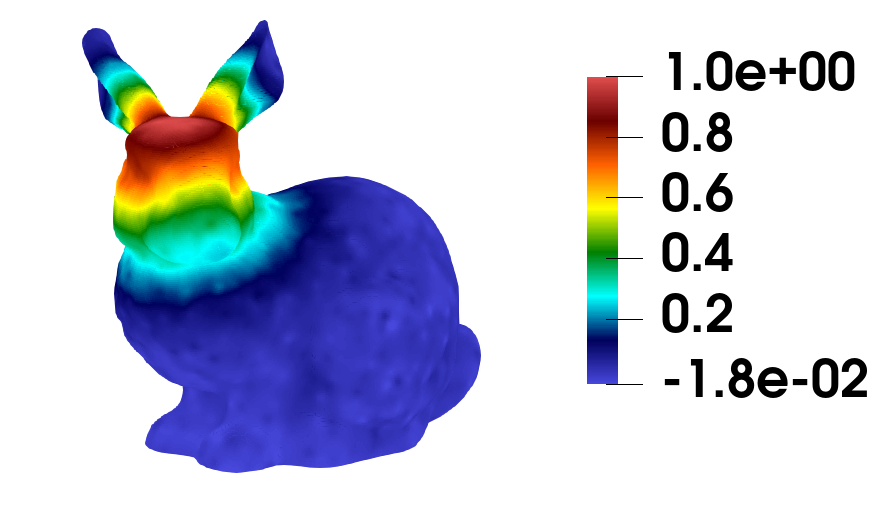} &
\includegraphics[width=0.3\linewidth]{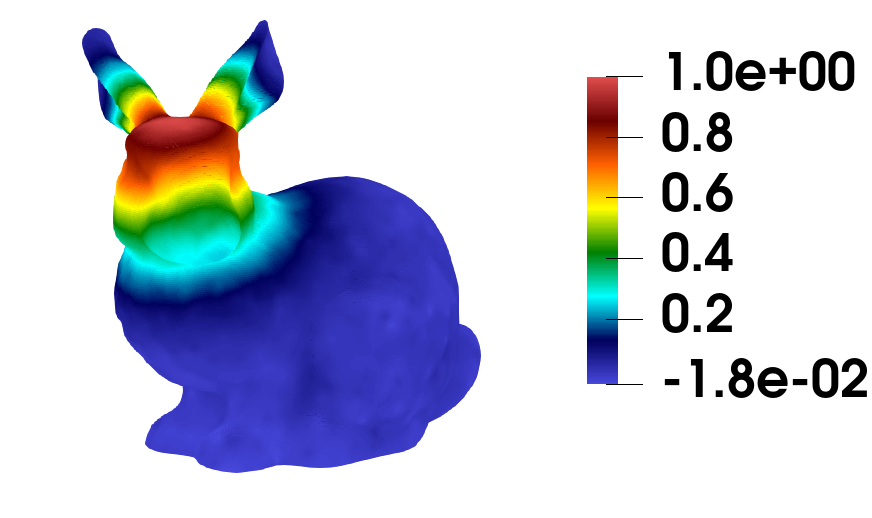} &
\includegraphics[width=0.3\linewidth]{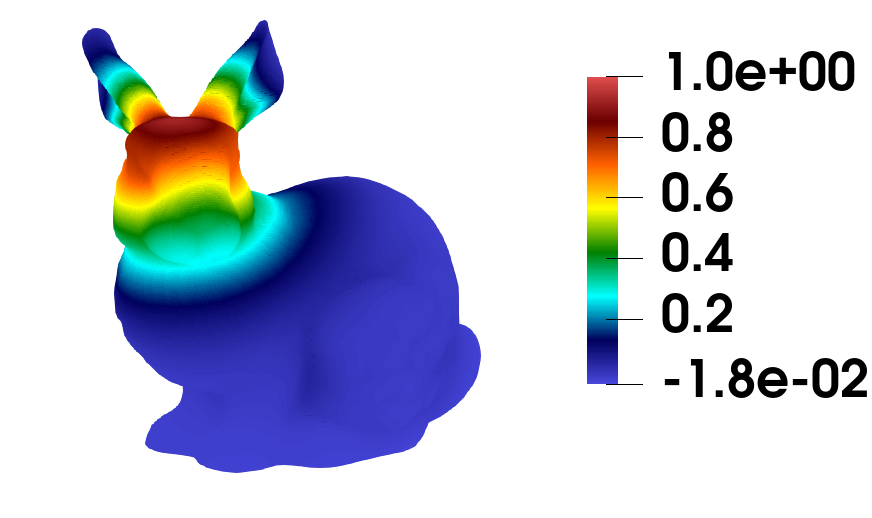}\\
\includegraphics[width=0.3\linewidth]{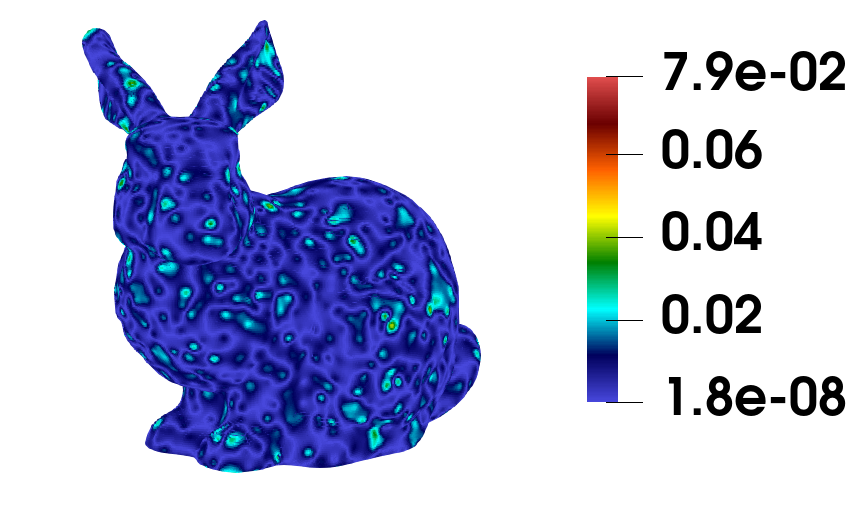} &
\includegraphics[width=0.3\linewidth]{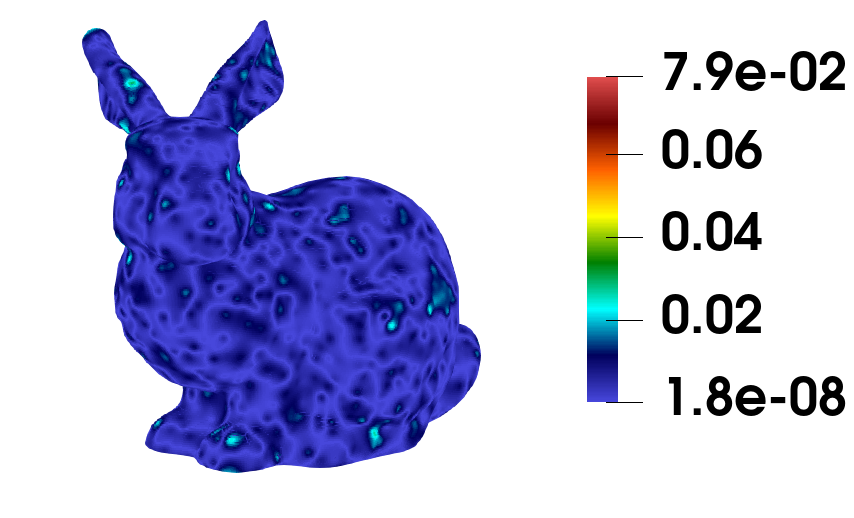} &
\includegraphics[width=0.3\linewidth]{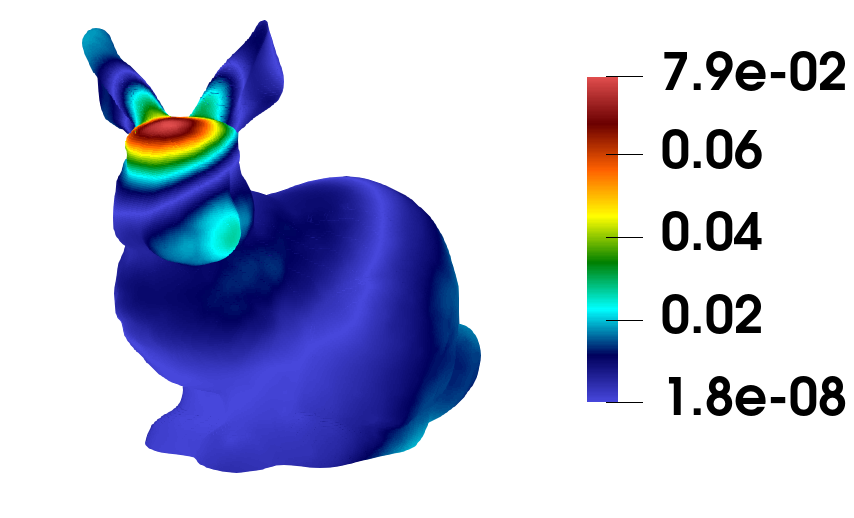}
\end{tabular}
\caption{Noisy boundary data on the Stanford bunny at the finest level
$\ell=6$, for $\Lambda_2=0$ (collocation), $\Lambda_2=10^{-2}$ and
$\Lambda_2=1$ (from left to right). Top row: approximant. Bottom row:
absolute pointwise error $\abs{u-u_\Phi}$, on the same color scale in the three
panels.}
\label{fig:bunny_pikl}
\end{figure}

\Cref{tab:bunny_pikl} collects the results and \Cref{fig:bunny_pikl} shows
the approximants and the errors at the finest level. Collocation interpolates
the noisy boundary values. With $\Lambda_2=10^{-2}$, the fit no longer
interpolates them: the error is below the collocation one on every level, and
the norm of the coefficient vector is less than half of the collocation one.
With $\Lambda_2=1$, the penalty dominates: the solution is smooth but biased,
and the error is larger than the collocation one on every level. 

\begin{table}[htb]
\centering\footnotesize
\begin{tabular}{rr ccc ccc ccc}
\toprule
 & & \multicolumn{3}{c}{$\Lambda_2=0$}
   & \multicolumn{3}{c}{$\Lambda_2=10^{-2}$}
   & \multicolumn{3}{c}{$\Lambda_2=1$}\\
\cmidrule(lr){3-5}\cmidrule(lr){6-8}\cmidrule(lr){9-11}
$\ell$ & $N_{\ell}$
 & $\varepsilon_2$ & $\norm{\bs\alpha}_2$ & it
 & $\varepsilon_2$ & $\norm{\bs\alpha}_2$ & it
 & $\varepsilon_2$ & $\norm{\bs\alpha}_2$ & it\\
\midrule
1 &      130 & $5.79\cdot10^{-2}$ & 3.18 & 2 & $5.14\cdot10^{-2}$ & 1.65 & 7 & $9.59\cdot10^{-2}$ & 0.97 & 15\\
2 &      454 & $2.91\cdot10^{-2}$ & 4.47 & 2 & $2.72\cdot10^{-2}$ & 1.79 & 6 & $3.38\cdot10^{-2}$ & 0.98 & 11\\
3 &   1\,737 & $1.81\cdot10^{-2}$ & 4.88 & 2 & $1.74\cdot10^{-2}$ & 1.92 & 4 & $3.16\cdot10^{-2}$ & 0.99 &  7\\
4 &   6\,745 & $1.66\cdot10^{-2}$ & 5.08 & 2 & $1.40\cdot10^{-2}$ & 2.01 & 3 & $3.38\cdot10^{-2}$ & 0.99 &  6\\
5 &  26\,469 & $1.55\cdot10^{-2}$ & 5.26 & 2 & $1.18\cdot10^{-2}$ & 2.10 & 2 & $3.49\cdot10^{-2}$ & 0.99 &  5\\
6 & 104\,589 & $1.46\cdot10^{-2}$ & 5.39 & 2 & $1.11\cdot10^{-2}$ & 2.17 & 2 & $3.48\cdot10^{-2}$ & 0.99 &  4\\
\bottomrule
\end{tabular}
\caption{PIKL on the Stanford bunny with five percent noise on the boundary
datum, for three values of the regularization weight: relative $\ell^2$ error
\eqref{eq:eps2}, norm of the coefficient vector and CG iterations of
\Cref{alg:pikl}.}
\label{tab:bunny_pikl}
\end{table}

\section{Conclusion}\label{sec:conclusion}
In this article, we have combined samplet compression with multiscale kernel
collocation for elliptic boundary value problems. In the samplet basis, the level
matrices have only $\mathcal{O}(N\log N)$ relevant entries, while the multiscale rescaling of the kernel controls the conditioning of the
rescaled interior blocks. Since the differential operator scales the interior
rows differently from the boundary rows, we have introduced a stable solver based
on a Schur-complement reduction, which is invariant under this scaling. We have
also derived a relation between the fill distances of the interior and boundary
points, together with the corresponding counting rule.

We have then considered physics-informed kernel learning, shown that it reduces
exactly to collocation as the regularization vanishes, and solved its normal
equations by a deflated conjugate gradient method that reuses the collocation
solver.

We have tested the method on Poisson problems in two and three dimensions, with
exact and noisy boundary data. In two dimensions, the compressed multiscale scheme
matches the accuracy of the dense one and reaches problem sizes far beyond those
of the dense scheme, while single scale collocation fails as the number of points
grows. On the Stanford bunny, the method remains efficient on a complex
three-dimensional geometry. In both cases, the number of conjugate gradient
iterations of the multiscale scheme stays bounded over the levels. For
physics-informed kernel learning, the experiments on exact data confirm the limit
to collocation, and the experiments with noisy boundary data show that a small
regularization reduces the error below that of collocation.
\bibliographystyle{plain}
\bibliography{bibliography}
\end{document}